\magnification =1098 \hsize 16truecm
\vsize 23truecm
\hfuzz =5pt
\scrollmode
\def\qed{\hbox{\vrule height 7pt depth 0pt width 7pt}}
\def\cqfd{\hfill\penalty 500\kern 10pt\qed\medbreak}
\def \e{{\varepsilon}}

\def \a{{\alpha}}

\def \qq{{\qquad}}
\def \b{{\beta}}
\def \s{{\sigma}}

\def \R{{\bf R}}
\def \Z{{\bf Z}}

\def \E{{\bf E}}

 \def \N{{\bf N}}

\def\be#1\ee{\begin{equation}#1\end{equation}}
\font\iit =cmmi7 at 7pt
\hfuzz =  5pt

\def \noi{{\noindent}}
  \font\ph=cmcsc10  at  10 pt
  \font\phh=cmcsc10  at  8 pt
  \font\iit =cmmi7 at 7pt
   \vskip 10mm
  \font\gum= cmbx10 at 11 pt
  \font\subgum= cmbx10 at 9 pt
  \font\ggum= cmbx10 at 10,7 pt
  \font\gem= cmbx10 at 10 pt

\def\ddate {\ifcase\month\or January\or
February\or March\or April\or May\or June\or July\or
August\or September\or October\or November\or December\fi\ {\the\day},\
{\sevenrm\the\year}}
 ${}$

\centerline{ {\ggum
Sampling the Lindel\"of Hypothesis with the Cauchy Random Walk}}

\vskip 1cm
 \centerline{Mikhail {\ph Lifshits}\ and\  Michel {\ph Weber}}
\vskip 0,6cm

{
\leftskip =1,9cm \rightskip=1,9cm  \noindent {\sl Abstract:\it We
study the behavior of the Riemann zeta function $\zeta({1\over 2} + it)$,
when $t$ is sampled by the Cauchy random walk. More precisely, let
$X_1,X_2,\ldots$ denote an infinite  sequence of independent Cauchy
distributed
random variables. Consider the sequence of partial sums
$S_n= X_1+ \ldots +X_n $, $n=1,2,\ldots$.  We   investigate the almost sure
asymptotic
behavior of the system
\vskip -5pt
$$  \zeta({1\over 2} + iS_n),  \qq n=1,2,\ldots$$
\vskip -3pt

\noi
We develop a complete second order theory for this system and show, by
using a classical
approximation formula of $\zeta(\cdot)$, that it behaves almost like a
system of
non-correlated variables. Exploiting this fact in relation with known criteria
for almost sure convergence, allows to prove the following almost sure
asymptotic behavior:
for any real $b>2$,
\vskip -5pt

$$\sum_{k=1}^n \zeta({1\over 2} + iS_k)
  \buildrel{(a.s.)}\over=
n +{\cal O}\big( n^{{1/2}}(\log n)^b\big)
$$
}
\par}

\footnote {}{{\iit Date}\sevenrm : on
\ddate\par
\vskip -2pt {\iit AMS}
\ {\iit Subject}\ {\iit Classification} 2000: Primary 11M06, 60G50,
Secondary 60F15.\par \vskip -2pt
  {\iit Keywords}:  Riemann zeta function, Lindel\"of Hypothesis, Cauchy
random walk. \par }
\rm

\vskip 0,7cm\noi \centerline{\gum
1. Introduction and Main Result}
\medskip\medskip

\noi Our work is devoted to the study of the celebrated Lindel\"of Hypothesis,
and our main theorems provide new quantitative results about the behavior
of the Riemann
zeta function along the critical line
$\Re s={1\over 2}$. Here and elsewhere we use the standard notation
$s=\sigma+it$ for the
complex argument. As is well-known, the Riemann zeta function defined
on the half-plane $\{s: \Re s>1\}$ by the series
$$\zeta(s)= \sum_{n=1}^\infty n^{-s}\leqno(1.1) $$
admits a meromorphic continuation to the entire complex plane, with the unique
and simple pole  of residue $1$ at $s=1$. In the
half-plane $\{s:\Re s\le 0\}$, the Riemann zeta function has simple zeros at
$-2,-4,-6,\ldots$, and only at these points which are called trivial zeros.
There exist also non-trivial zeros in the band  $\{s: 0< \Re s < 1\}$.
We refer for these basic facts for instance to [Bl] (Propositions
IV.10 \& IV.11, p.84).

Two great conjectures are related to the behavior of $\zeta(s)$. The
Riemann Hypothesis
(RH) asserts that all non-trivial zeros of the function $\zeta$ have
abscissa ${1\over 2}$, while Lindel\"of Hypothesis (LH) claims that
$$\zeta( {1\over 2}+it) ={\cal O}(t^\e) \leqno(1.2)
$$
for every positive $\e$; or, what turns out to be equivalent ([T], Chap. XIII
p.276), that
$$\zeta( \s+it) ={\cal O}(t^\e) \leqno(1.3)
$$
for every positive $\e$ and every $\s \ge {1\over 2}$. The validity of
RH  implies ([T], Theorem 14.14, Chap. XIV p.300) that
$$\zeta( {1\over 2}+it) ={\cal O}\left(\exp\left\{ A{\log t\over \log\log
t}\right\}\right), \leqno(1.4)
$$
$A$ being a constant, which is even a stronger form of
LH; the latter  being strictly weaker than RH.

There are various equivalent reformulations of the LH. Here we follow [T]
Chap. XIII,
and recall that the validity of (1.3) is equivalent to any of the three
following assertions
$${1\over T}\int_1^T \big|\zeta({1\over 2} + it)\big|^{2k} {\rm d}t= {\cal
O}\big( T^\e\big), \quad\qq k=1,2,\ldots \leqno(1.5)
$$
$${1\over T}\int_1^T \big|\zeta(\s + it)\big|^{2k} {\rm d}t= {\cal O}\big(
T^\e\big), \quad\qq \s>{1\over 2}, \quad  k=1,2,\ldots \leqno(1.6)
$$
$$\lim_{T\to\infty}
{1\over T}\int_1^T \big|\zeta(\s + it)\big|^{2k} {\rm d}t
=
\sum_{n=1}^\infty {d_k^2(n)\over n^{2\s}}, \quad\qq \s>{1\over 2},
\quad  k=1,2,\ldots
\leqno(1.7)
$$
where $d_k(n)$ denotes the number of representations of integer $n$ as a
product of
$k$ factors. There are some classical results related to (1.7).
For every positive integer $k>2$, it is known ([T], Theorem 7.7 p.125)
that
$\lim_{T\to
\infty}{1\over T}\int_1^T
\big|\zeta(\s + it)\big|^{2k} {\rm d}t=
\sum_{n=1}^\infty {d_k^2(n)\over n^{2\s}}$ if $\s>1-1/k $.
The same result also holds ([T], Theorem 7.11 p.132) for
non-integer $k$ such that $0<k\le 2$ and $\s>{1/2}$,
and is proved by using a theorem of Carlson.
\medskip

to briefly
the various

The study of the LH has been over the last
century and up to now, the object of continuous and considerable
efforts of numerous mathematicians, not exclusively number
theorists, but also of probabilists,
starting from important contributions of P\'olya [Po]. Up to now, the
best known result towards (1.3) is due to Huxley [H2]
     $$\zeta( {1\over 2}+it) ={\cal O}(t^{32/205 +\e}), \qq (\forall
\e>0)\leqno(1.8)
$$
  and $32/205=0,156097561..$. Regarding the equivalent formulation in
terms of power moments (1.5), there is the following
  satisfactory estimate (see [I1] Theorem 5.1 p.129) due to Ingham for the
case $k=2$:
$$ \int_1^T \big|\zeta({1\over 2} + it)\big|^{4} {\rm d}t=
(2\pi^2)^{-1}T\log^4T + {\cal O}(T\log^3 T)   . \leqno(1.9)
$$
Beyond this case, for instance for $k=3$, nothing comparable has been
proved yet, and we may just cite the following much weaker estimate
$$ \int_1^T
\big|\zeta({1\over 2} + it)\big|^{2k} {\rm d}t= {\cal O}\big(
T^{(k+2)/4}\log ^{C(k)} T\big),
\quad\qq 2\le k\le 6, \leqno(1.10)
$$
 $C(k)$ is a constant depending on $k$. Define the (modified) Mellin transform
of the zeta function:
$${\cal M}_k(s) = \int_1^\infty |\zeta({1\over 2}+iu)|^{2k}u^{-s} {\rm d}u,
\quad\qq k\in \N, \s=\Re s\ge c(k)>1, \leqno(1.11) $$
 where $c(k)$ is such a constant for which the integral in (1.9) converges
absolutely.   It has been recently proved by Ivi\'c
([I2], Corollary 1) that the validity of LH is also equivalent to the
property that for every $k\in \N$,  ${\cal M}_k(s)$ is regular
for $\s>1$ and satisfies ${\cal M}_k(1+\e+it)\ll_{k,\e}1 $.

The LH has also a connection with the function $S(T)$ where we recall ([T],
Section 9.3) that $S(T)$
denotes the value of
$$ \pi^{-1}\arg\zeta({1\over 2} + iT) , \qq ( \arg\zeta(s)=
\arctan{\Im\zeta(s)\over \Re\zeta(s)} )
$$
obtained by continuous variation along the straight lines joining $2$,
$2+iT$, ${1\over 2}+iT$, starting with the value $0$. Whereas it
is known that
$S(T)={\cal O}(\log T)$, the validity of LH would imply $S(T)=o(\log T)$,
see ([T]   Theorem 9.4 p.181 and Theorem 13.6 p.281)
respectively. In [Gh], Ghosh answering a question raised by Selberg,
studied the value distribution of the modulus $|S(t)|$ and showed
that

$${\rm meas}\Big\{T < t < T+H:  |S(t)|<\sigma \sqrt{\log\log t}   \Big\}
          ={1\over \sqrt{2\pi}}  \int_{-\sigma}^\sigma e^{-x^2} d x + o(1)
. H, \leqno(1.12) $$
is valid for $ T^\alpha < H < T$ and any fixed $\alpha > 1/2$.  He also proves
that on the Riemann Hypothesis this result holds for any fixed $\alpha >
0$.

There are numerous results focused on the value distribution of the zeta
function, since the seminal work of Bohr and Jessen [BJ].
We refer for instance to the book  of Joyner [Jo].
The limiting value distribution (so called Bohr-Jessen measure) has been
extensively studied in the works of Hattori and Matsumoto [HM],
as well as by Laurin\v cikas, who proved in [La] (see also his joint work with
Steuding [LS]) by means of probabilistic methods that the LH is actually
equivalent to the fact that for arbitrary positive reals
$\e$ and $a$
$${1\over T} {\rm meas}\bigg\{ t\in [0,T] :\big|\zeta({1\over
2}+it)\big|<xT^\e\bigg\}=1-{\cal O}\Big({\Delta(T)\over
1+x^a}\Big)\leqno(1.13)$$  holds for all $x$ large enough, where $\Delta(T)$
is an arbitrary function such that $\Delta(T)=o(1)$.

We refer to [T], [I] and the recent survey [GM] for other related
results, like for instance ([T] Chap. XIII) the
relationships between the LH and the distribution of the zeros of the
zeta function. It is worth noticing that the validity of the RH has other
interesting consequences  concerning the
asymptotic behavior of $\zeta(1+it)$.
For instance, Theorem 14.8 p.290 of [T] implies that
$$ \big|\log \zeta (1+it)\big|\le \log\log\log t +A,\leqno(1.14)
$$
$A$ being a constant,
 whereas Vinogradov [V] proved $\zeta(1+it) ={\cal O}((\log t)^{2/3})$
\smallskip

Now we would like to mention some probabilistic methods involved in
the study of the zeta function. In [Bi], [BPY], various identities in
distribution linking functionals of Brownian motion with the elliptic
theta function
$$\Theta(u) =\sum_{n\in \Z} e^{-\pi n^2u}, \leqno(1.15)
$$
allow to reinterpret or retrieve differently the famous functional
equation (see [H1], Chap. 11, Eq. (11.3) and (11.7) or [T] Chap. II
or else [Bl] Part. 5, Chap.3 p.136) valid for any complex $s$
$$\eqalign{  \pi^{-{1\over2}s}\Gamma({1\over2}s)\zeta (s)  &  = \pi^{-{
1\over2}(1-s)}\Gamma({1-s\over2} )\zeta (1-s)\cr     &  =\int_1^\infty
{1\over2} \big(\Theta(x)-1\big)\big(
x^{{1\over2}s-1}+x^{-{1\over2}s-{1\over2}}\big)\, {\rm d}x-\big\{
s(1-s)\big\}^{-1}  .\cr}\leqno(1.16)
$$
In [W1], another form of the functional equation, linking the zeta function
with the value distribution of the divisors of the spin random walk, has been
recently established. The approach of [W1] is different
from those of the above quoted papers.
\smallskip

To conclude this brief description of involved probabilistic methods, it
seems necessary to mention the, although not relevant in the
present work, very actively developing random matrix theory, modelling
the pair correlation of zeros of the zeta function.
We refer for instance to [Bi] (Section 2) for a short
glimpse to this theory based on the striking observation made by Dyson
that the asymptotic
distribution formula for the distances between the zeros of the
zeta function proposed by Montgomery [Mo] exactly describes
the distribution of the distances between the eigenvalues of a Gaussian
random Hermitian matrix (recall that
Hilbert and P\'olya suggested that the zeros of the
zeta function are likely the eigenvalues of some Hilbertian
self-adjoint operator). This is also motivated by the analogy existing
between the explicit formula for the zeros given by  A. Weil, and Selberg
trace formula for the discrete eigenvalues of the Laplace
operator in the hyperbolic half-plane ($\Delta= -y^2\Big( { \partial^2\over
\partial x^2}+{ \partial^2\over \partial y^2}\Big)$).
None of the reviewed approaches will be, however, implemented in the
present work.

Here, our aim is to study the asymptotic behavior of the zeta function
along the critical line $\s={1\over 2}$ by modelling the
time $t$  with the Cauchy random walk. Let $X_1,X_2,\ldots$ denote an
infinite
sequence of independent Cauchy distributed random variables
(with characteristic function $\varphi(t) =e^{-|t|}$), then the time $t$ is
modelled by the
sequence of partial sums
$$S_n= X_1+ \ldots +X_n .\leqno(1.17)
$$
In order to understand the behavior of $\zeta({1\over 2} + it)$ when $t$ tends
to infinity, we propose to investigate the almost sure asymptotic behavior
of the system
$$  \zeta_n:=\zeta({1\over 2} + iS_n), \quad\qq n=1,2, \ldots \leqno(1.18)
$$
Put for any positive integer $n$
$${\cal Z}_n = {\zeta}({1/  2}+iS_n) -\E\,  {\zeta}({1/ 2}+iS_n)
=\zeta_n-\E \zeta_n. \leqno(1.19)
$$
We develop a complete second order theory of the system $\{ {\cal Z}_n,n\ge
1\}$.
The main striking fact we obtain is that the this system nearly behaves
like a system of non-correlated
variables, i.e. the variables ${\cal Z}_n$ are weakly orthogonal. More
precisely, we prove
 \medskip

\noi {\gem Theorem 1.} {\it There exist constants $C,C_0$ such that}
$$\eqalign{\E\,   | {\cal Z}_n|^2
&=\log n + C
+ o(1), \qquad n\to\infty,
\cr
\hbox{\it and for $m>n+1$, \qq}\big|\E\,{\cal Z}_n\overline{{\cal Z}_m}\big|
&\le C_0 \max \Big( {1\over n},{1\over 2^{ m-n  }}   \Big). \qq\qq\qq\qq
\cr}
$$
\medskip


\noi  {\gem Remark.}  The explicit value of $C$ is
$$
  C=C_E -2  + 2\int_0^1 \phi(\alpha) d\alpha +
    2 \int_1^\infty \left( \phi(\alpha) -   {1\over 2\alpha} \right)
    d\alpha,
$$
where $C_E$ is the Euler constant and
$\phi(\alpha)={\alpha e^{\alpha} -2 e^{\alpha}+\alpha+2 \over  2\alpha^2
(e^{\alpha}-1)}$.
\smallskip

Exploiting Theorem 1 in the context of the
known criteria  for almost sure convergence, we prove
the following theorem, which displays a rather slow growth of the
zeta function on the critical line, when sampled by the Cauchy random walk.
\medskip

 \noi {\gem Theorem 2.} {\it For any real $b>2$,
$$\lim_{n\to \infty}{\sum_{k=1}^n \zeta({1\over 2} + iS_k) - n  \over
n^{{1/2}}(\log n)^b
}\buildrel{(a.s.)}\over=0 ,$$
and}
$$ \bigg\|\sup_{n\ge 1}{\big|\sum_{k=1}^n \zeta({1\over 2} + iS_k) - n
\big|\over n^{{1/2}}(\log
n)^b }\bigg\|_2<\infty .
$$
The used notation $a.s.$ (for almost surely) means that the corresponding
property holds with probability one.
  \medskip

\noi  {\gem Remark.}  We believe that the results similar to Theorems 1 and
2 are valid
for sampling with a large class of random walks with discrete or continuous
steps. Quite surprisingly, the necessary moment expressions we obtain for
Cauchy distribution
are by far more explicit (which made our project feasible) than in other
cases,
e.g. for Gaussian or Bernoulli distributions.
 \medskip

Our approach is based on the following classical approximation result
(see for instance Theorem 4.11 p.67 in   [T]): letting, as usually,
$s=\s+it$, we have
$$ \zeta(s) =\sum_{n\le x} {1\over n^s} -{x^{1-s}\over 1-s} +{\cal O}(
x^{-\s}),\leqno(1.20) $$
uniformly for $\s\ge \s_0>0$, $|t|\le T_x:=2\pi x/C$, $C$ is any constant
$>1$.

\noi  Therefore, the second order theory of the system $({\cal Z}_n)$
follows from
a study of the same kind concerning the auxiliary system
$$ Z_n(x) =  Z_n = \sum_{k\le x} {1\over k^{\sigma+iS_n}}
-    {x^{1-(\sigma+iS_n)}\over 1-(\sigma+iS_n)},   \quad\qq n=1,2\ldots, \
x>0.
\leqno(1.21)
$$
\medskip

\noi  The investigation of $Z_n(x)$ occupies the whole Section 2, and
constitutes
the main part of the technical work. In the concluding Section 3 we show
that $Z_n(x)$
approximates zeta function well enough and prove Theorems 1 and 2.
\smallskip

\medskip\medskip\medskip\medskip

\centerline{\gum 2. Second order theory of (Z\lower +1,2mm\hbox{\subgum n})}

\bigskip \noi
We begin with some basic notation. We write $Z_n(x)=Z_n=Z_{n1}-Z_{n2}$
with
$$
Z_{n1} =Z_{n1}(x)= \sum_{k\le x} {e^{-i(\log k) S_n}\over k^\sigma},
\leqno(2.1)$$
\vskip -17,8 pt
$$
 Z _{n2}= Z_{n2}(x)=  {x^{1-\sigma} e^{-i(\log x) S_n}\over
1-(\sigma+iS_n)} . \leqno(2.2)$$
In order to investigate the covariance structure, we study the behavior of
the first and second order moments of
$ Z_n $, and  the correlation $\E\,  Z_n\bar Z_m $, from which are easily
derived the second order distances $\E\, |Z_n-Z_m|^2$,
$m>n$.
We write
$$\eqalign{\E\, |Z_{n}|^2&= E|Z_{n1}|^2+ E|Z_{n2}|^2 -2\, \Re\, \E\,
Z_{n1}\bar Z_{n2}
\cr
\E\,  Z_n\bar Z_m
&=\E\,  Z_{n1} \bar Z_{m1} -\E\, Z_{n1} \bar Z_{m2} -\E\,  Z_{n2} \bar
Z_{m1} +\E\,  Z_{n2} \bar Z_{m2},\cr}\leqno(2.3)
$$
The following integral representation will be used repeatedly
$$
  {1\over 1-s}=\int_0^1 u^{-s} {\rm d}u,  \qq \qq \Re s< 1.
$$
For the first moments, we have for $x\ge 1$,
$$\eqalign{
\E\, Z _{n2}
&=    x^{1-\sigma}\int_0^1   \E\, e^{-i(\log x) S_n}  e^{-(\log
u)(\sigma+iS_n)}du
 =    x^{1-\sigma}\int_0^1 u^{-\s}  \E\, e^{-i(\log xu) S_n}   du
\cr
&=    x^{1-\sigma}\int_0^1 {e^{-|\log(xu)|n}du \over u^\s}
=x^{1-\sigma}\left( \int_0^{1/x} {(xu)^n du \over u^\s}
               +   \int_{1/x}^1 {(xu)^{-n}du \over u^\s}  \right)
\cr
&= x^{1-\sigma}\left( {x^{\s-1}\over n-\s+1} +{x^{\s-1}-x^{-n}\over n+\s-1}
\right)
= {2n\over n^2-(1-\s)^2} -{x^{1-\s-n}\over n+\s-1}\ ,
\cr}
$$
and
$$\E\, Z_{n1} =  \sum_{k\le x}{ 1\over k^\sigma} \E\, e^{-i(\log k) S_n}
=\sum_{k\le x}{ 1\over k^{\sigma+n}}.
$$
Therefore,
$$
   \E\, Z _{n }=\E\, Z _{n1}-\E\, Z _{n2}
   =\sum_{k\le x}{ 1\over k^{\sigma+n}} -  {2n\over n^2-(1-\s)^2} +
{x^{1-\s-n}\over n+\s-1}
\buildrel{x\rightarrow\infty}\over\longrightarrow \zeta(\s+n) -   {2n\over
n^2-(1-\s)^2}  ,\leqno(2.4)
$$
for any integer $n$ and $\s>0$.
\smallskip  In subsequent calculations we will find explicit and asymptotic
formulas
for   $\E\,|Z_{n1}|^2$ (see (2.21)),
$\E\,|Z_{n2}|^2$ (see (2.12)), $\E\, Z_{n2}\bar Z_{m2}$ (see (2.13)), $\E\,
Z_{n1}\bar Z_{m2}$  (see (2.15)), $\E\, Z_{m1}\bar Z_{n2}$
(see (2.17)),
$\E\, Z_{n1}\bar Z_{n2}$ (see (2.16)).
The final answers are given in Section 2.5.
\medskip

\noi{\gum 2.1.  Exact formulae related to $\bf Z_{n2}$}
\medskip
\noi We begin with proving three exact
formulae stated in the following proposition.
\medskip

\noi{\gem Proposition 1.} {\it For $m=n$ and for $m>n+1$ we have
 $$ \E\, Z_{n2}\bar Z_{m2} = A+ B x^{-n+(1-\sigma)} +C
x^{-(m-n)+2(1-\sigma)}, \leqno(2.5)
$$
where
$$\eqalign{
A&=   {4n(m-n)\over ((m-n)^2-4(1-\sigma)^2)(n^2-(1-\sigma)^2)},
\cr B &=   {2(m-n)\over (2n-m+(1-\sigma))(m+(1-\sigma))(n-(1-\sigma))},
\cr C &=   {3n-m+2(1-\sigma)\over
(2n-m+(1-\sigma))(2(1-\sigma)-(m-n))(n+(1-\sigma))}.
\cr}$$
For all $m\ge n$ we have
$$\displaylines{
 (2.6)\qq\E\, Z_{n1}\bar Z_{m2} = \sum_{k\le x}
\Bigl[
{-2(m-n)k^{-n-\sigma}\over (m+(1-\sigma))(2n-m+(1-\sigma))}
\hfill\cr\hfill +
{2nk^{-(m-n)+1-2\sigma}\over(m-(1-\sigma))(2n-m+(1-\sigma))}
-
{k^{n-\sigma}x^{-m+(1-\sigma)}\over m-(1-\sigma)}
\Bigr],
\cr}$$
and}
$$\E\,  Z_{m1}\bar Z_{n2} =
\sum_{k\le x} \left(
{2n k^{-(m-n)+1-2\sigma}\over n^2-(1-\sigma)^2}
-
{k^{2n-m-\sigma}x^{-n+(1-\sigma)}\over n-(1-\sigma)}
\right).\leqno(2.7)
$$
\medskip\noi{\it Proof.}
We start with the proof of (2.5). Recall that
$$ Z_{n2} =
{x^{1-\sigma}e^{-i\log x S_n}\over 1-(\sigma + i S_n)},
\qq\qq \bar Z_{m2} =
{x^{1-\sigma}e^{i\log x S_m}\over 1-(\sigma - i S_m)}.
$$
Thus,
$$
\E\,  Z_{n2}\bar Z_{m2}
= x^{2(1-\sigma)} \E \left[
  e^{i\log x (S_m-S_n)}
\ {1\over 1-(\sigma + i S_n)}
\ {1\over 1-(\sigma - i S_m)}
                     \right].
$$
Using again the integral representation
 $
  {1\over 1-s}=\int_0^1 u^{-s} {\rm d}u$, $s\ae 1$,
 we obtain
$$\eqalign{\E\, \Big[   e^{i\log x (S_m-S_n)}
\  {1\over 1-(\sigma + i S_n)}& \ {1\over 1-(\sigma - i S_m)}
     \Big]
\cr & = \int_0^1\int_0^1 u^{-\sigma}v^{-\sigma}
\E\, e^{i(\log x+\log v)(S_m-S_n)+i(\log v- \log u)S_n} {\rm d}u {\rm d}v\cr &=
 \int_0^1\int_0^1 u^{-\sigma}v^{-\sigma}
 e^{-|\log x+\log v|(m-n)-|\log v- \log u|n} {\rm d}u {\rm d}v.
\cr}$$
Next, we split the square $[0,1]^2$ in four domains.
\smallskip

\noi -- For the first domain, $u\le v, 1/x\le v$, we have
$$\eqalign{  \int_{1/x}^1 {\rm d}v \int_0^v u^{-\sigma}v^{-\sigma}
 & e^{-|\log x+\log v|(m-n)-|\log v- \log u|n} {\rm d}u
\cr &   =
\int_{1/x}^1 {\rm d}v \int_0^v u^{-\sigma}v^{-\sigma}
 x^{-(m-n)} v^{-(m-n)}(u/v)^n {\rm d}u
\cr & =
x^{-(m-n)} \int_{1/x}^1 v^{-m-\sigma} {\rm d}v
    \int_0^v  u^{n-\sigma}{\rm d}u
  =
 {x^{-(m-n)}\over n+(1-\sigma)} \int_{1/x}^1 v^{-m-\sigma+n+(1-\sigma)}
{\rm d}v
\cr & = {x^{-(m-n)}\over n+(1-\sigma)} \ \cdot \
 {x^{(m-n)- 2(1-\sigma)}-1\over (m-n) - 2(1-\sigma)}
\cr & =
 {x^{-2(1-\sigma)}- x^{-(m-n)}\over ((m-n)-2(1-\sigma))(n+(1-\sigma))}.
\cr}$$
Remark that this calculation does not go through in the case
$m=n+1,\sigma=1/2$ that we excluded. The same is valid for many other
subsequent formulas but we will not stress this fact anymore.

Thus, for the first domain,
$$\displaylines{  (2.8)\qq x^{2(1-\sigma)} \int_{1/x}^1 {\rm d}v \int_0^v
u^{-\sigma}v^{-\sigma}
 e^{-|\log x+\log v|(m-n)-|\log v- \log u|n} {\rm d}u
\hfill\cr\hfill
= {1 - x^{-(m-n)+2(1-\sigma)}\over ((m-n)-2(1-\sigma))(n+(1-\sigma))}.
\cr}$$

\noi -- For the second domain, $u \le v \le 1/x$, we have
$$\eqalign{ \int_0^{1/x}{\rm d}v \int_0^v u^{-\sigma}v^{-\sigma}
& e^{-|\log x+\log v|(m-n)-|\log v- \log u|n} {\rm d}u
\cr &=
\int_0^{1/x} {\rm d}v \int_0^v u^{-\sigma}v^{-\sigma}
 x^{m-n} v^{m-n}(u/v)^n {\rm d}u
\cr &=
 x^{m-n} \int_{0}^{1/x} v^{m-2n- \sigma} {\rm d}v
    \int_0^v  u^{n-\sigma} {\rm d}u
\cr &=
  {x^{m-n}\over n+(1-\sigma)} \int_0^{1/x} v^{m-2n-\sigma+n+(1-\sigma)}
{\rm d}v
\cr &=
  {x^{m-n}\over n+(1-\sigma)} \ \cdot \
 {x^{-(m-n)- 2(1-\sigma)}\over (m-n)+2(1-\sigma)}
\cr &=
  {x^{-2(1-\sigma)}\over ((m-n)+2(1-\sigma))(n+(1-\sigma))}.
\cr}$$

Thus, for the second domain,
$$\displaylines{(2.9)\qq x^{2(1-\sigma)} \int_0^{1/x} {\rm d}v \int_0^v
u^{-\sigma}v^{-\sigma}
 e^{-|\log x+\log v|(m-n)-|\log v- \log u|n} {\rm d}u
\hfill\cr\hfill=
 {1\over ((m-n)+2(1-\sigma))(n+(1-\sigma))}.
\cr}$$

\noi -- For the third domain, $u \ge v \ge 1/x$, we have
$$\eqalign{ \int_{1/x}^1 {\rm d}v \int_v^1 u^{-\sigma}v^{-\sigma}
& e^{-|\log x+\log v|(m-n)-|\log v- \log u|n} {\rm d}u
\cr & =
\int_{1/x}^{1} {\rm d}v \int_v^1 u^{-\sigma} v^{-\sigma}
 (xv)^{-(m-n)} (v/u)^n {\rm d}u
\cr & =
 x^{-(m-n)} \int_{1/x}^{1} v^{2n-m- \sigma}{\rm d}v
    \int_v^1  u^{-n-\sigma} {\rm d}u
\cr & =
   {x^{-(m-n)}\over n-(1-\sigma)}
  \int_{1/x}^{1} v^{2n-m-\sigma}(v^{-n+(1-\sigma)}-1) {\rm d}v
\cr & =
 {x^{-(m-n)}\over n-(1-\sigma)} \
\left(
 {x^{(m-n)- 2(1-\sigma)}\over (m-n)-2(1-\sigma)}
+
 {x^{-(2n-m+(1-\sigma))}\over 2n-m+(1-\sigma)}
\right)
\cr & \quad
 -
  {x^{-(m-n)}\over n-(1-\sigma)}
\left(
 {1\over (m-n)-2(1-\sigma)}
+
 {1\over 2n-m+(1-\sigma)}
\right).
\cr}$$

Thus, for the third domain,
$$\eqalign{
  x^{2(1-\sigma)} \int_{1/x}^{1} {\rm d}v \int_v^1  u^{-\sigma}v^{-\sigma}
& e^{-|\log x+\log v|(m-n)-|\log v- \log u|n} {\rm d}u
\cr &=
 {1\over ((m-n)-2(1-\sigma))(n-(1-\sigma))}
\cr &=
-  {x^{-n+(1-\sigma)}\over (2n-m+(1-\sigma))(n-(1-\sigma))}
\cr &=
+ \
 {x^{-(m-n)+2(1-\sigma)}\over (2n-m+(1-\sigma))((m-n)-2(1-\sigma))}.
\cr}\leqno(2.10)$$

\noi -- For the fourth and the last domain, $u \ge v, 1/x\ge v$, we have
$$\eqalign{
 \int_0^{1/x} {\rm d}v \int_v^1 u^{-\sigma}v^{-\sigma}
& e^{-|\log x+\log v|(m-n)-|\log v- \log u|n} {\rm d}u
\cr &=
\int_0^{1/x}  {\rm d}v \int_v^1 u^{-\sigma} v^{-\sigma}
 (xv)^{(m-n)} (v/u)^n {\rm d}u
\cr &=
x^{(m-n)} \int_0^{1/x} v^{m- \sigma}{\rm d}v
    \int_v^1  u^{-n-\sigma} {\rm d}u
\cr &=
   {x^{(m-n)}\over n-(1-\sigma)}
  \int_0^{1/x} (v^{m-n+1-2\sigma}-v^{m-\sigma}) {\rm d}v
\cr &=
 {x^{(m-n)}\over n-(1-\sigma)} \
\left(
 {x^{-(m-n)-2(1-\sigma)}\over (m-n)+2(1-\sigma)}
-
 {x^{-m-(1-\sigma))}\over m+(1-\sigma)}
\right).
\cr}$$

Thus, for the fourth domain,
$$\displaylines{(2.11)\qq
  x^{2(1-\sigma)} \int_0^{1/x} {\rm d}v \int_v^1  u^{-\sigma}v^{-\sigma}
 e^{-|\log x+\log v|(m-n)-|\log v- \log u|n} {\rm d}u
\hfill\cr \hfill= {1\over (n-(1-\sigma))((m-n)+2(1-\sigma))}
 -\  {x^{-n+(1-\sigma)}\over (n-(1-\sigma))(m+(1-\sigma))}.
\cr}$$
By summing up eight terms in (2.8),
 (2.9), (2.10), (2.11),
we arrive at  (2.5).

\medskip An important particular case of (2.5) is $m=n$,
where

$$\E\,Z_{n2}\bar Z_{n2} =    {x^{2(1-\sigma)}\over (1-\sigma)(n+(1-\sigma))}.
\leqno(2.12)$$\bigskip

Now we pass to the proof of (2.6). By the definition,
$$\eqalign{\E\, Z_{n1}\bar Z_{m2} &= \E\, \sum_{k\le x} k^{-\sigma}
e^{-i\log k S_n}
{x^{1-\sigma}e^{i\log x S_m}\over 1-(\sigma-iS_m)}
\cr &=
x^{1-\sigma}\ \sum_{k\le x} k^{-\sigma} \E\,  \left[
  e^{-i\log k S_n + i\log x S_m}
\int_0^1 v^{-\sigma+iS_m} {\rm d}v \right]
\cr &=
 x^{1-\sigma}\ \sum_{k\le x} k^{-\sigma}
     \int_0^1 v^{-\sigma} E  e^{-i\log k S_n+i(\log x +\log v)S_m} {\rm d}v
\cr &=
x^{1-\sigma}\ \sum_{k\le x} k^{-\sigma}
     \int_0^1 v^{-\sigma}
              e^{-|\log(xv)| (m-n) -|\log(xv/k)|n} {\rm d}v.
\cr}$$
We calculate the last integral by splitting $[0,1]$ in three intervals.
First,
$$\eqalign{
 \int_0^{1/x} v^{-\sigma}
  e^{-|\log(xv)|(m-n) -|\log(xv/k)|n} {\rm d}v
  &=
\int_0^{1/x} v^{-\sigma} (xv)^{(m-n)} (xv/k)^n {\rm d}v
\cr &=
{x^m\over k^n} \int_0^{1/x} v^{m-\sigma} {\rm d}v
\cr &=
{x^m\over k^n}\ \cdot \ { x^{-m-(1-\sigma)}\over m+(1-\sigma)}
.
\cr}$$
Second,
$$\eqalign{ \int_{1/x}^{k/x} v^{-\sigma}
  e^{-|\log(xv)| (m-n) -|\log(xv/k)|n} {\rm d}v
  &=
\int_{1/x}^{k/x} v^{-\sigma} (xv)^{-(m-n)} (xv/k)^n {\rm d}v
\cr &=
{x^{2n-m}\over k^n}  \int_{1/x}^{k/x} v^{2n-m-\sigma} {\rm d}v
\cr &=
{x^{2n-m}\over k^n}\ \cdot \
{ x^{-2n+m-(1-\sigma)} (k^{2n-m+(1-\sigma)}-1) \over 2n-m+(1-\sigma)}
\cr &=
{x^{-(1-\sigma)}\over 2n-m+(1-\sigma)}\ \cdot \
\left( k^{-(m-n)+(1-\sigma)}- k^{-n} \right).
\cr}$$

Third,
$$\eqalign{  \int_{k/x}^{1} v^{-\sigma}
  e^{-|\log(xv)| (m-n) -|\log(xv/k)|n} {\rm d}v
&=
\int_{k/x}^{1} v^{-\sigma} (xv)^{-(m-n)} (k/xv)^n {\rm d}v
\cr &= {k^n\over x^m}  \int_{k/x}^{1} v^{-m-\sigma} {\rm d}v
\cr &=  {k^n\over x^m(m-(1-\sigma))} \left( (k/x)^{-m+(1-\sigma)}-1 \right)
\cr &=  {k^{-(m-n)+(1-\sigma)} \over x^{1-\sigma}(m-(1-\sigma))}
-  {k^n\over x^m(m-(1-\sigma))}.
\cr}$$
By summing up three answers, multiplying by $k^{-\sigma}$, adding up over $k$,
and multiplying by $x^{1-\sigma}$, we easily arrive at (2.6).
\medskip

The proof of (2.7) is very similar. Indeed, we have
by the definition,
$$\eqalign{\E\, Z_{m1}\bar Z_{n2}  &=\E\,\sum_{k\le x} k^{-\sigma}
e^{-i\log k S_m}
{x^{1-\sigma}e^{i\log x S_n}\over 1-(\sigma-iS_n)}
\cr &= x^{1-\sigma}\ \sum_{k\le x} k^{-\sigma} \E\, \left[
  e^{-i\log k S_m+i\log x S_n}
\int_0^1 v^{-\sigma+iS_n} {\rm d}v \right]
\cr &= x^{1-\sigma}\ \sum_{k\le x} k^{-\sigma}
     \int_0^1 v^{-\sigma} \E\, e^{-i\log k S_m+i(\log x +\log v)S_n} {\rm d}v
\cr &= x^{1-\sigma}\ \sum_{k\le x} k^{-(m-n)-\sigma}
     \int_0^1 v^{-\sigma} e^{-|\log(xv/k)|n} {\rm d}v.
\cr}$$
We calculate this integral by splitting $[0,1]$ in two intervals.
First,
$$\eqalign{ \int_0^{k/x} v^{-\sigma} e^{-|\log(xv/k)|n} {\rm d}v
  &= \int_0^{k/x} v^{-\sigma} (xv/k)^{n} {\rm d}v
\cr &= {x^n\over k^n}\ \cdot \ {(k/x)^{n+(1-\sigma)}\over n+(1-\sigma)}
\cr &= {k^{(1-\sigma)}\over x^{(1-\sigma)}(n+(1-\sigma))}.
\cr}$$
Second,
$$\eqalign{  \int_{k/x}^1 v^{-\sigma} e^{-|\log(xv/k)|n} {\rm d}v
  &= \int_{k/x}^1 v^{-\sigma} (k/xv)^{n} {\rm d}v
\cr &= {k^n\over x^n}\, \cdot \,
 {1\over n-(1-\sigma)}
\left(  {x^{n-(1-\sigma)}\over k^{n-(1-\sigma)}} -1 \right)
\cr &= {k^{(1-\sigma)}\over x^{(1-\sigma)}(n-(1-\sigma))}
-
{k^{n}\over x^{n}(n-(1-\sigma))}.
\cr}$$
By summing up two answers, multiplying by $k^{-(m-n)-\sigma}$,
adding up over $k$, and multiplying by $x^{1-\sigma}$,
we easily arrive at (2.7).

\medskip\noi {\gum 2.2. Asymptotic formulae related to $\bf Z_{n2}$}
\medskip\noi
Here we give a brief asymptotic analysis of the results
obtained in previous section regarding the behaviour of
exact expressions at $x\to\infty$. For the sake of brevity,
we only consider $\sigma=1/2$.
\medskip
    It follows immediately from (2.5) that for $m>n+1$
$$\E\, Z_{n2}\bar Z_{m2} =
 {4n(m-n)\over ((m-n)^2-1)(n^2-1/4)}+o(1), \qquad x\to \infty,
\leqno(2.13)$$
while (2.12) yields
$$
\E\, Z_{n2}\bar Z_{n2} =
 {2x\over n+1/2}.
\leqno(2.14)$$

\noi Next, (2.6) implies
$$\eqalign{
\E\, Z_{n1}\bar Z_{m2} &=
{-2(m-n)\zeta(n+1/2)\over (m+1/2)(2n-m+1/2)}
+ o(1)
\cr &\quad + \
{2n\over (m-1/2)(2n-m+1/2)}
\sum_{k\le x} k^{-(m-n)}
  - \
 {x^{-m+1/2}\over m-1/2}
\sum_{k\le x}  k^{n-1/2}.
\cr}$$
Remark that for $m>n+1$ the second term converges and the second
one is negligible, since
$$
x^{-m+1/2} \sum_{k\le x}k^{n-1/2}
\le x^{-m+1/2} \cdot x \cdot x^{n-1/2} = x^{-(m-n)+1}= o(1).
$$
Hence, for $m>n+1$, we obtain
$$
\E\, Z_{n1}\bar Z_{m2}=
{-2(m-n)\zeta(n+1/2)\over (m+1/2)(2n-m+1/2)}
 + \
{2n \ \zeta(m-n)\over (m-1/2)(2n-m+1/2)}
+ o(1), \qquad x\to \infty.
\leqno(2.15)$$
  When $m=n>2$, we use that, by second order Euler--Maclaurin formula,
$$
 \sum_{k\le x} k^{n-1/2} =
    {x^{n+1/2}\over n+1/2} +  {x^{n-1/2}\over 2} +
  o\left( x^{n-1/2} \right)
$$
and obtain
$$\eqalign{
\E\,  Z_{n1}\bar Z_{n2} &=
{2nx\over n^2-1/4} - {x\over n^2-1/4}-{1\over 2n-1}
+ o(1)
\cr &=  {2x\over n+1/2}-{1\over 2n-1}
+ o(1), \qq\qquad x\to \infty.
\cr}\leqno(2.16)$$
Now let us consider (2.7) that now writes
$$
\E\, Z_{m1}\bar Z_{n2} =
{2n\over n^2-1/4} \sum_{k\le x} k^{-(m-n)}
-
{x^{-n+1/2}\over n-1/2}
\sum_{k\le x} k^{2n-m-1/2} .
$$
When $m>n+1$, the first term converges and the second one is
vanishing, since
$$
x^{-n+1/2} \sum_{k\le x} k^{2n-m-1/2}
\le x^{-n+1/2} \cdot x \cdot x^{2n-m-1/2} = x^{-(m-n)+1}= o(1).
$$
Thus, we get
$$\E\, Z_{m1}\bar Z_{n2} = {2n \ \zeta(m-n)\over n^2-1/4}
+ o(1), \qq \qquad x\to \infty.
\leqno(2.17)$$
 On the other hand, putting $m=n $ in (2.7), yields again
(2.16).
\medskip


\noi{\gum 2.3. Calculation of $\E\,\bf Z_{n1}\bar Z_{m1}, m>n+1$}
\medskip\noi
Let us fix $\sigma\in[1/2,1)$ and $m,n$ so that $m\ge n$.
We have
$$\eqalign{
\E\, Z_{n1}\bar Z_{m1}&=  \E\, \sum_{k,l\le x}
{e^{-i\log k S_n}\over k^\sigma}\
{e^{i\log l S_m}\over l^\sigma}
\cr & =
 \E\, \sum_{k,l\le x}   {1\over k^\sigma l^\sigma}
e^{i(\log l -\log k) S_n}\
e^{i\log l (S_m-S_n)}
  =
\sum_{k,l\le x}   {1\over k^\sigma l^\sigma}
\left( {\min(k,l)\over \max(k,l)}\right)^n \ l^{-(m-n)}
\cr & =  S_1+S_2+S_0,
\cr}\leqno(2.18)$$
where
$$
\eqalign{ S_1 & = \sum_{k\le x} k^{n-\sigma} \sum_{l=k+1}^x l^{-m-\sigma},
\cr S_2 &= \sum_{l\le x} l^{2n-m-\sigma} \sum_{k=l+1}^x k^{-n-\sigma}
    = \sum_{k\le x} k^{2n-m-\sigma} \sum_{l=k+1}^x l^{-n-\sigma},
\cr S_0 &= \sum_{k\le x} k^{-(m-n)-2\sigma}.
\cr}$$
We specify this to $\sigma=1/2$, so that
$$\eqalign{
S_1 & = \sum_{k\le x} k^{n-1/2} \sum_{l=k+1}^x l^{-m-1/2},
\cr S_2 &= \sum_{k\le x} k^{2n-m-1/2} \sum_{l=k+1}^x l^{-n-1/2},
\cr
S_0 &= \sum_{k\le x} k^{-(m-n)-1}.
\cr}\leqno(2.19)$$
For $m>n+1$ we obviously have
$$
S_0 =  \zeta((m-n)+1) + o(1), \qquad x \to \infty.
$$
Next,
$$
S_1 = \sum_{k=1}^\infty  k^{n-1/2} \sum_{l=k+1}^\infty l^{-m-1/2}
      +o(1),   \qquad x \to \infty.
$$
Moreover, for $m-n>1$,
$$\eqalign{
 \sum_{k=1}^\infty  k^{n-1/2} \sum_{l=k+1}^\infty l^{-m-1/2}
  & =
\sum_{k=1}^\infty  k^{n-1/2}  \ \theta_{k,m}
\int_{k}^\infty u^{-m-1/2}{\rm d}u
\cr & =
 {\theta_{m,n} \over m-1/2}
\sum_{k=1}^\infty  k^{-(m-n)}
\cr & =
 { \theta_{m,n} \over m-1/2}
\ \zeta(m-n).
\cr}$$
 Here and elsewhere $\theta$'s are different constants in $[0,1]$.
\smallskip
  Exactly in the same way we obtain
$$
S_2 = \sum_{k=1}^\infty  k^{2n-m-1/2} \sum_{l=k+1}^\infty l^{-n-1/2}
      +o(1),   \qquad x \to \infty,
$$
and
$$
\sum_{k=1}^\infty  k^{2n-m-1/2} \sum_{l=k+1}^\infty l^{-n-1/2}
= {\theta'_{m,n}\over n-1/2} \ \zeta(m-n).
$$
Thus, finally, for $m>n+1$
$$\eqalign{  \E\, Z_{n1}\bar Z_{m1}
 &= \zeta((m-n)+1)+\theta \left(  {1\over m-1/2} +  {1\over n-1/2}
\right)\zeta(m-n) +o(1), \qq x \to \infty,
\cr}\leqno(2.20)$$
with $\theta=\theta(n,m)\in [0,1]$.

\medskip\noi {\gum 2.4. Calculation of $\E\, \bf  Z_{n1}\bar Z_{n1}$.}
\medskip\noi
Our aim is to prove the   following
formula
$$\E\, Z_{n1}\bar Z_{n1}= {2x\over n+1/2} +K_n + o(1),
\qquad x\to\infty,
\leqno(2.21)$$
with
$$K_n= \log n + C + o(1), \quad n\to \infty
$$
and
$$
C=C_E -1+
2\int_0^1 \phi(\alpha) d\alpha +
2 \int_1^\infty \left( \phi(\alpha) -   {1\over 2\alpha} \right)
d\alpha,
\leqno(2.22)$$
where $C_E$ is the Euler constant and
$\phi(\alpha)={\alpha e^{\alpha} -2 e^{\alpha}+\alpha+2 \over  2\alpha^2
(e^{\alpha}-1)}$.
\smallskip

Let us start proving (2.21). We already know
that
$$\eqalign{
\E\, Z_{n1}\bar Z_{n1}=
\sum_{k,l\le x}  {1\over k^{1/2}l^{1/2}}
\left( {\min(k,l)\over \max(k,l)}\right)^n
 =  2 \sum_{l\le x}  {1\over l^{n+1/2}}
      \sum_{k\le l} k^{n-1/2}
  - \sum_{l\le x}  {1\over l}.
\cr}$$
We use Euler-Maclaurin formula of the first order:
$$
 \sum_{k\le l} k^{n-1/2} =
{l^{n+1/2}-1\over n+1/2}
+ {l^{n-1/2}+1\over 2}
+ \sum_{k\le l-1} A_k,
$$
where
$$ A_k= (n-3/2)\int_0^1 (k+t)^{n-3/2}(t-1/2) {\rm d}t.
$$
By summing up we arrive at
$$\eqalign{
 2 \sum_{l\le x}  {1\over l^{n+1/2}}
      \sum_{k\le l} k^{n-1/2}
  & =
2 \sum_{l\le x}  {1\over l^{n+1/2}}
\left(
{l^{n+1/2}-1\over n+1/2}
+ {l^{n-1/2}+1\over 2}
+ \sum_{k\le l-1} A_k
\right)
\cr & =
  {2x\over n+1/2} +  \sum_{l\le x} {1\over l} +
2 \left({1\over 2}-{1\over n+1/2} \right)
 \sum_{l\le x} {1\over l^{n+1/2}}
\cr & \quad +\
2 \sum_{l\le x}  {1\over l^{n+1/2}}
      \sum_{k\le l-1} A_k
\cr & =
  {2x\over n+1/2} +  \sum_{l\le x} {1\over l} +
{n-3/2\over n+1/2}\ \zeta(n+1/2)
\cr & \quad
+
2 \sum_{k=1}^\infty A_k
 \sum_{l=k+1}^{\infty}  {1\over l^{n+1/2}}
+o(1), \qq\qquad x\to \infty.
\cr}$$
Hence,
$$\eqalign{
\E\, Z_{n1}\bar Z_{n1} =
  {2x\over n+1/2} + {n-3/2\over n+1/2}\ \zeta(n+1/2)
 +
2 \sum_{k=1}^\infty A_k
 \sum_{l=k+1}^{\infty}  {1\over l^{n+1/2}}
+o(1), \qq x\to \infty.
\cr}\leqno(2.23)$$
Now, it remains to analyze the behavior of the double sum
$$
 S= \sum_{k=1}^\infty A_k
 \sum_{l=k+1}^{\infty}  {1\over l^{n+1/2}}
$$
when $n\to\infty$. Let denote
$$\eqalign{
B_k&=B_k(n)=\int_0^1 (k+t)^{n-3/2}(t-1/2) {\rm d}t,
\cr
D_k&=D_k(n)=  \sum_{l=k+1}^{\infty}  {1\over l^{n+1/2}},\cr
D'_k & =D'_k(n)=  \sum_{l=k+2}^{\infty}  {1\over l^{n+1/2}}.
\cr}$$
Then we have
$$\eqalign{
S &= (n-3/2) \sum_{k=1}^\infty B_k D_k
\cr & =
(n-3/2) \bigg[
  \sum_{k=n}^\infty B_k D_k  +
  \sum_{k=1}^{n-1} B_k \left( D'_k +(k+1)^{-n-1/2}\right)\bigg].
\cr}$$
We will show in the next section that
$$\lim_{n\to\infty}
(n-3/2) \sum_{k=n}^\infty B_k D_k
=\int_0^1 \phi(\alpha) d\alpha,
\leqno(2.24)$$
where $\phi(\alpha)={\alpha e^{\alpha} -2 e^{\alpha}+\alpha+2 \over
2\alpha^2 (e^{\alpha}-1)}$;
$$
\lim_{n\to\infty}
(n-3/2) \sum_{k=1}^{n-1} B_k D'_k
=\int_1^\infty \phi_1(\alpha) d\alpha,
\leqno(2.25)$$
where $\phi_1(\alpha)=
  {\alpha -2 +\alpha e^{-\alpha} +2 e^{-\alpha}\over  2\alpha^2
(e^{\alpha}-1)}$; and
$$
 \lim_{n\to\infty}
\bigg( (n-3/2) \sum_{k=1}^{n-1} B_k (k+1)^{-n-1/2}
-\sum_{k=1}^{n-1}  {1\over 2(k+1)}
\bigg)
=\int_1^\infty \phi_2(\alpha) d\alpha,
\leqno(2.26)$$
where $\phi_2(\alpha)=
       {2 e^{-\alpha}+\alpha e^{-\alpha}-2\over 2\alpha^2}$.
Note that
$$\phi_1(\alpha)+\phi_2(\alpha)=
       {2+ 2\alpha - 2 e^{\alpha}\over 2\alpha^2 (e^{\alpha}-1)}
=\phi(\alpha)-  {1\over 2\alpha}\ .
$$
It follows from (2.24), (2.25), (2.26) that
$$ S= {1\over 2} \sum_{k=1}^{n-1}  {1\over k+1}
+ \int_0^1 \phi(\alpha) d\alpha +
\int_1^\infty \left( \phi(\alpha) -   {1\over 2\alpha} \right)
d\alpha +o(1).
$$
Recall that
$$
\sum_{k=1}^{n-1}  {1\over k+1} =
\sum_{k=1}^{n}  {1\over k} - 1 =
\log n +C_E - 1 + O(1/n)
$$
where $C_E$ is the Euler constant.
Thus we finally obtain
$$ 2S = \log n +C_E - 1 + O(1/n)
+ 2\int_0^1 \phi(\alpha) d\alpha +
2 \int_1^\infty \left( \phi(\alpha) -   {1\over 2\alpha} \right)
d\alpha + o(1),
$$
as asserted in (2.21).


\medskip
\noi{\gum 2.5. Proofs of (2.24), (2.25), (2.26).}
\medskip\noi We want to show (2.24), i.e.
$$\lim_{n\to\infty}
 n  \sum_{k=n}^\infty B_k D_k
=\int_0^1 \phi(\alpha) d\alpha,
$$
with $\phi(\alpha)={\alpha e^{\alpha} -2 e^{\alpha}+\alpha+2 \over
2\alpha^2 (e^{\alpha}-1)}$.
\smallskip
To achieve this, we obtain that for any (large) fixed $M>1$, uniformly over
$k\in [n,Mn]$, it is true that
$$ nB_kD_k\sim \int_{n/(k+1)}^{n/k} \phi(\a)d\a.
\leqno(2.27)$$
Since $\phi$ is uniformly continuous, we have
$$ \int_{n/(k+1)}^{n/k} \phi(\a)d\a \sim \phi({n\over k+1})({n\over k
}-{n\over k+1})\sim \phi(\b_k){n\over (k+1)^2},$$
where $\b_k={n\over k+1}\in ]0,1]$. Thus we need to check
$$B_kD_k\sim \phi(\b_k){1\over (k+1)^2},
$$
or, equivalently
$$(k+1)^2 B_kD_k\sim  \phi(\b_k).
\leqno(2.28)$$
We have
$$(k+1)^2 B_kD_k= nB_k(k+1)^{-(n-{1\over 2})} \ \cdot \ D_k(k+1)^{
(n+{1\over 2})} \ \cdot \ {k+1\over n},
\leqno(2.29)$$
and will show that
$$   nB_k(k+1)^{-(n-{1\over 2})} \sim    {1\over 2}(1+e^{-\b_k})+
{e^{-\b_k}-1\over \b_k},
\leqno(2.30)$$
$$  D_k(k+1)^{ (n+{1\over 2})}  \sim   (1-e^{-\b_k})^{-1}.\leqno(2.31)$$
Since
$$\big( {1\over 2}(1+ e^{-\b})+   {e^{-\b }-1\over \b }\big)
\big(1-e^{-\b}\big)^{-1}
\ \cdot \
\b^{-1}=\big(\b e^\b+\b+2-2e^\b\big)\big(2\b^2(e^\b-1)\big)^{-1}=\phi(\b),
$$
(2.28) would follow from (2.29)-(2.31). Now we prove (2.30). We have, by
variable change $t=1-{k+1\over n}v$,

$$\eqalign{ B_k& =  \int_0^1(k+t)^{n-3/2}(t-{1\over 2})dt\cr &
 ={(k+1) \over n}(k+1)^{n-3/2}\int_0^{n/(k+1)}(1-v/n)^{n-3/2}\big({1\over
2}-{k+1\over n}
v\big) dv
\cr &\sim {(k+1)^{n-1/2}\over n}\int_0^{\b_k}e^{-v}\big({1\over 2}-{v\over
\b_k}\big) dv .\cr}$$
By using the explicit formula
$$ \int_0^\b e^{-v}({1\over 2}-{v\over \b })dv={1\over 2}(1+e^{-\b})+
{e^{-\b}-1\over 2},$$
we arrive at (2.30).
\smallskip
Now we prove (2.31). We have
$$D_k(k+1)^{n+1/2}=\sum_{h=1}^\infty \big({k+1\over
k+h}\big)^{n+1/2}=\sum_{h=1}^\infty \big(1+ {h-1\over k+1}\big)^{-(n+1/2
)}= \sum_{h=0}^\infty \big(1+ {h \over k+1}\big)^{-(n+1/2 )}.
\leqno(2.32)$$
By using (2.32), we have
$$\eqalign{  D_k(k+1)^{n+1/2}&\sim \sum_{h=0}^\infty e^{-h(n+1/2)/(k+1)} \cr
&\sim \big(1-\exp(- {n+1/2\over k+1})\big)^{-1}\sim \big(1-\exp(- {n \over
k+1})\big)^{-1}\cr &\sim
 \big(1-\exp(- \b_k)\big)^{-1},\cr}$$
as asserted in (2.31).
\smallskip Now (2.27) is proved completely and we obtain, for any fixed $M$
$$\liminf_{n\to \infty} n\sum_{k=n}^\infty B_kD_k\ge \lim_{n\to \infty}
\sum_{k=n}^{Mn}\int_{n/(k+1)}^{n/k}\phi(\a)d\a=\int_{1/M}^{1}\phi(\a)d\a.$$
By sending $M$ to infinity, we arrive at
$$\liminf_{n\to \infty} n\sum_{k=n}^\infty B_kD_k\ge \int_{0}^{1}\phi(\a)d\a.
\leqno(2.33)$$
Similarly, we get for any $M>1$
$$\limsup_{n\to \infty} n\sum_{k=n}^{Mn} B_kD_k\le
\int_{1/M}^{1}\phi(\a)d\a\le \int_{0}^{1}\phi(\a)d\a.
\leqno(2.34)$$
Thus we only need to show that
$$\lim_{M\to \infty} \ \limsup_{n\to \infty} n\sum_{k>Mn} B_kD_k=0.
\leqno(2.35)$$
Then (2.34) and (2.35) will imply
$$\limsup_{n\to \infty} n\sum_{k=n}^\infty B_kD_k\le
\int_{0}^{1}\phi(\a)d\a, $$
and thus finish the proof of (2.24), being coupled with (2.33)
\smallskip
We now prove (2.35). We use again  that
$$\eqalign{ B_k& ={k+1\over n}(k+1)^{n-3/2}\int_0^{{n\over k+1}}(1-{v\over
n})^{n-3/2}\big( {1\over 2}-{k+1\over n}v\big) dv
\cr & ={(k+1)^{n-1/2}\over n} \int_0^{\b_k}(1-{v\over n})^{n-3/2}\big(
{1\over 2}-{v\over \b_k} \big) dv,\cr}$$
and observe that
$$\eqalign{
\Big|\int_0^\b (1-{v\over n})^{n-3/2}\big( {1\over 2}-{v\over \b  } \big)
dv\Big|
&=
\Big|\int_0^\b \Big((1-{v\over
n})^{n-3/2}-1\Big) \  \big( {1\over 2}-{v\over \b  } \big)  dv
+ \int_0^\b \big( {1\over 2}-{v\over \b  } \big) dv\Big|
\cr
& =\Big|\int_0^\b \Big((1-{v\over n})^{n-3/2}-1\Big) \big( {1\over
2}-{v\over \b  } \big)dv+\ 0\Big|
\cr
&\le  \int_0^\b \Big|(1-{v\over
n})^{n-3/2}-1\Big| dv \cr}.
$$
As
$$ \Big|(1-{v\over
n})^{n-3/2}-1\Big|= (n-3/2)\int_{1- {v\over n}}^1 y^{n-5/2}dy\le
(n-3/2){v\over n}\le v, $$
we get
$$ B_k\le {(k+1)^{n-1/2}\over n}\b_k^2.
\leqno(2.36)$$
 Similarly, we will prove that
$$ D_k.(k+1)^{n+1/2}\le C\b_k^{-1}.
\leqno(2.37)$$
It follows that
$$nB_kD_k\le (k+1)^{n-1/2}\b_k^2(k+1)^{-(n+1/2)}
\ \cdot \
C\b_k^{-1}= C(k+1)^{-1}\b_k=C{n\over (k+1)^2},
$$
whence
$$n\sum_{k>Mn}B_kD_k\le Cn\sum_{k>Mn}{1\over (k+1)^2}\le Cn/(Mn)\le C/M,
$$
and (2.35) follows. Thus it remains to check (2.37). Recall that by (2.32)
$$D_k \ (k+1)^{n+1/2}=\sum_{h=0}^\infty\big(1+{h\over k+1}\big)^{-(n+1/2)}. $$
We split the sum in two: firstly, by using
$$1+s\ge e^{s\log 2}, \qq\qq 0\le s\le 1,
$$
we have
$$\eqalign{\sum_{h=0}^{k+1} \big(1+{h\over k+1}\big)^{-(n+1/2)}&\le
\sum_{h=0}^\infty\exp\big(-(n+1/2){h\over k+1}\log
2\big)=\Big(1-\exp\big(- {n+1/2\over k+1}\log 2\big)\Big)^{-1}
\cr& \le \Big(1-\exp\big(- 4\b_k\big)\Big)^{-1}\le C\b_k^{-1}, \cr}$$
for all $0\le \b_k\le 1$.

Secondly,
$$\eqalign{\sum_{h>k+1}\big(1+{h\over k+1}\big)^{-(n+1/2)}{(k+1)\over
(k+1)}&\le (k+1)\int_1^\infty (1+x)^{-(n+1/2)}dx={(k+1)\over
(n-1/2)}2^{-(n-1/2)}\cr& \le 2^{3/2 -n}{(k+1)\over
n}\le C\b_k^{-1},\cr}$$
and we are done with (2.37) and with all the proof of (2.24).
\bigskip

The proof of (2.25) is completely similar to that of (2.24). We want to
show that
$$\lim_{n\to\infty}
(n-3/2) \sum_{k=1}^{n-1} B_k D'_k
=\int_1^\infty \phi_1(\alpha) d\alpha,
$$
with $\phi_1(\alpha)=
  {\alpha -2 +\alpha e^{-\alpha} +2 e^{-\alpha}\over  2\alpha^2
(e^{\alpha}-1)}$.
\smallskip

 The main point is that for any (large) fixed $M>1$, uniformly over $k\in
[{n\over M}, n]$, we have
$$ nB_k D'_k\sim \int_{n/(k+1)}^{n/k} \phi_1(\a)d\a.
$$
By continuity of $\phi_1$, we have
$$\int_{n/(k+1)}^{n/k} \phi_1(\a)d\a\sim \phi_1({n\over k+1})\big( {n\over
k}- {n\over k+1}\big)\sim \phi_1(\b_k){n\over (k+1)^2},$$
where $\b_k={n\over k+1}\in [1,M]$. Thus we need to check
$$ B_kD'_k\sim \phi_1(\b_k)/(k+1)^2,
$$
or, equivalently,
$$ (k+1)^2 B_kD'_k\sim \phi_1(\b_k)  .
\leqno(2.38)$$
We write
$$(k+1)^2 B_kD'_k= nB_k(k+1)^{-(n-1/2)} \ \cdot \ D'_k(k+1)^{n+1/2}\ \cdot
\ {k+1\over n}.
\leqno(2.39)$$
Next, we use again (2.30) which claims
$$ nB_k(k+1)^{-(n-1/2)}\sim {1\over 2}(1+e^{-\b_k})+ {e^{-\b_k}-1\over \b}.
\leqno(2.40)$$
Moreover, we obtain from (2.31) that
$$\eqalign{D'_k(k+1)^{n+1/2}&= \Big( D_k-(k+1)^{-(n+1/2)}\Big)
(k+1)^{n+1/2}=D_k(k+1)^{n+1/2}-1\cr& \sim
\big(1-e^{-\b_k})^{-1}=e^{-\b_k}/\big(1-e^{-\b_k}\big).\cr}
\leqno(2.41)$$
Since
$$\Big[{1\over 2}\big(1+e^{-\b}\big)+
{e^{-\b}-1\over\b}\Big]\big({e^{-\b}\over 1-e^{-\b}} \big){1\over
\b}=\Big[\b\big(1+e^{-\b}\big)+2( e^{-\b}-1)\Big]{1\over 2\b^2(e^\b
-1)}=\phi_1(\b ),$$
we obtain (2.38) from (2.39) via (2.40) and (2.41). We derive  next from
(2.38) that for any fixed $M>1$
$$ \liminf_{n\to\infty} n\sum_{k=n/M}^nB_kD'_k\ge \lim_{n\to\infty}
n\sum_{k=n/M}^n\int_{n/(k+1)}^{n/k}\phi_1(\a)d\a=\int_1^M
\phi_1(\a)d\a.$$
By sending $M$ to infinity, we arrive at
$$\liminf_{n\to\infty} n\sum_{k=1}^nB_kD'_k\ge \lim_{M\to \infty}
\liminf_{n\to\infty} n\sum_{k=n/M}^nB_kD'_k \ge \int_1^\infty
\phi_1(\a)d\a.
\leqno(2.42)$$
Similarly, we get for any $M>1$
$$\limsup_{n\to\infty} n\sum_{k=n/M}^nB_kD'_k\le \int_1^M
\phi_1(\a)d\a\le   \int_1^\infty
\phi_1(\a)d\a.
\leqno(2.43)$$
Thus the only thing we need to show is
$$ \lim_{M\to \infty} \limsup_{n\to\infty} n\sum_{k=1}^{ n/M} B_kD'_k =0.
\leqno(2.44)$$
Then (2.43) and (2.44) will imply
$$\limsup_{n\to\infty} n\sum_{k=1}^nB_kD'_k\le  \int_1^\infty
\phi_1(\a)d\a,$$
and this, after coupling with (2.42), will finish the proof of (2.25).
\smallskip
We still have, by (2.36)
$$ nB_k\le (k+1)^{n-1/2}\b_k^2,
\leqno(2.45)$$
and will now evaluate $D'_k$ as follows
$$\eqalign{ D'_k \ \cdot \ (k+1)^{n+1/2}&=\sum_{h=2}^\infty \big({k+1\over
k+h}\big)^{n+1/2}=\sum_{h=2}^\infty \big(1+{h-1\over
k+1}\big)^{-(n+1/2)} =\sum_{h=1}^\infty \big(1+{h \over k+1}\big)^{-(n+1/2)}
\cr &=\left(\sum_{h=1}^{k+1}+\sum_{h= k+2}^\infty\right) \big(1+{h \over
k+1}\big)^{-(n+1/2)}
.\cr}$$
By using again $1+s\ge e^{s\log 2}$, $0\le s\le 1$, we have
$$\eqalign{ \sum_{h=1}^{k+1}  \big(1+{h \over k+1}\big)^{-(n+1/2)}&\le
\sum_{h=1}^\infty \exp\Big\{-{h(n+1/2).\log 2\over
k+1}\Big\}\cr &\le
\exp\big\{-{n \over k+1} \ \cdot \ \log 2\big\}\Big(1-\exp \{-{n \over k+1}
\}\Big)^{-1}
\cr &\le C
\exp\big\{-{n \over k+1}\ \cdot \ \log 2\big\} \le C2^{-\b_k},\cr}$$
for all $k\le n$.
\smallskip
We also have
$$\eqalign{\sum_{h= k+2}^\infty \big(1+{h \over k+1}\big)^{-(n+1/2)}
{k+1\over k+1}&\le (k+1)\int_1^\infty (1+x)^{-(n+1/2)} dx
\cr & ={k+1\over n-1/2}2^{-(n-1/2)} \le 4 \, \cdot \, 2^{-n}\le
C2^{-n/(k+1)}=C2^{-\b_k}.\cr}$$
It follows that
$$D'_k(k+1)^{n+1/2}\le C 2^{-\b_k},$$
and by (2.45),
$$\eqalign{ nB_kD'_k&=nB_k(k+1)^{-(n+1/2)}D'_k (k+1)^{  n+1/2 }\le (k+1)^{
n+1/2 }\b_k^2(k+1)^{-(n+1/2)}
\ \cdot \ C \ \cdot \ 2^{-\b_k}\cr &
={\b_k^2\over (k+1)} \ \cdot \ C \ \cdot \ 2^{-\b_k}\le C{\b_k \over
(k+1)}2^{-\b_k/2}=C{n\over (k+1)^2}2^{-\b_k/2}.\cr}
$$
We finally obtain
$$n\sum_{k=1}^{n/M}B_kD'_k\le  C\sum_{k=1}^{n/M}{n\over
(k+1)^2}2^{-\b_k/2}\le C\int_M^\infty 2^{-x/2}dx\ \to \ 0,$$
as $M$ tends to infinity, as claimed in (2.44); so that (2.25) is proved
completely.
\bigskip

Finally, we prove (2.26). By definition,
$B_k=\int_0^1(k+t)^{n-3/2}(t-{1\over 2})dt$, and we have to investigate the
limit
behavior of the sum
$$\sum_{k=1}^{n-1}(k+1)^{-(n+1/2)}
B_k=\sum_{k=1}^{n-1}\int_0^1{(k+t)^{n-3/2}\over (k+1)^{n-3/2}(k+1)^2}
(t-{1\over 2})dt.$$
By the variable change $t=1-{k+1\over n}v$, we come to
$$\eqalign{ \sum_{k=1}^{n-1}{k+1\over n}&\int_0^{n/(k+1)} \Big[{
(k+1)-{(k+1)\over n}v\over  k+1 }\Big]^{n-3/2}{\big({1\over 2}-{ k+1
\over n}v\big)\over (k+1)^2}  dv\cr & =  \sum_{k=1}^{n-1} {1\over
2(k+1)n}\int_0^{n/(k+1)}\big(1-{v\over n}\big)^{n-3/2}dv - {1\over
n^2}\sum_{k=1}^{n-1} \int_0^{n/(k+1)}\big(1-{v\over n}\big)^{n-3/2}vdv.\cr}
 \leqno(2.46)$$
We show for the second term
$$ {(n-3/2)\over n^2}\sum_{k=1}^{n-1}
 \int_0^{n/(k+1)}\big(1-{v\over n}\big)^{n-3/2}vdv\ \buildrel{n\to
\infty}\over {\longrightarrow}\ \int_0^1e^{-v}vdv + \int_1^\infty
e^{-v} dv= 1-e^{-1}.
\leqno(2.47)$$
Note that by integration by parts
$$\int_0^1e^{-v}vdv=-\int_0^1 vd(e^{-v})=-\Big[ve^{-v}\big|_0^1
-\int_0^1e^{-v}dv \Big]=-\big[e^{-1}-(1-e^{-1})\big]=1-2e^{-1}.$$
Since $\int_1^\infty e^{-v} dv=e^{-1}$, (2.47) follows.
\smallskip
Write the sum from (2.47) as one integral:
$${1\over
n }\sum_{k=1}^{n-1}
 \int_0^{n/(k+1)}\big(1-{v\over n}\big)^{n-3/2}vdv={1\over n} \int_0^\infty
\big(1-{v\over n}\big)^{n-3/2}\#\big\{ k:k+1\le n,\,
k+1\le {n\over v}\big\}vdv,$$
then split the integral over the domains $[0,1]$ and $]1,\infty[$, getting
$$\int_0^1 \big(1-{v\over n}\big)^{n-3/2}{\#\big\{ k:k+1\le n \big\}\over
n}vdv+ \int_1^\infty \big(1-{v\over n}\big)^{n-3/2}{\#\big\{
k:  k+1\le {n\over v}\big\}\over n}vdv .$$
It is obvious that the first integral converges to $\int_0^1e^{-v} vdv$ and
the second one to $\int_1^\infty e^{-v} dv$, since in both
cases the theorem of dominated convergence applies. Therefore (2.47) is proved.
\smallskip

\noi Now consider the second expression in (2.46). After multiplying by
$(n-3/2)$ we get
$${n-3/2\over n}\bigg(\sum_{k=1}^{n-1}   {X_k\over 2(k+1)}-\sum_{k=1}^{n-1}
{Y_k\over 2(k+1)} +\sum_{k=1}^{n-1}   {Z_k\over
2(k+1)}\bigg),
\leqno(2.48)$$
where
$$X_k=\int_0^\infty e^{-v} dv =1, \quad Y_k=\int_{n/(k+1)}^\infty e^{-v}
dv=e^{-n/(k+1)},\quad Z_k=\int_0^{n/(k+1)} \Big(\big(1-{v\over
n}\big) -e^{-v}\big) dv. $$
Obviously, the first sum equals
$$\sum_{k=1}^{n-1}   {1\over 2(k+1)} +{\cal O}\big( {\log n\over
n}\big)=\sum_{k=1}^{n-1}   {1\over 2(k+1)} +o(1),$$
as asserted in (2.26). For the second term in (2.48), we have
$$\sum_{k=1}^{n-1}   {Y_k\over 2(k+1)}=\sum_{k=1}^{n-1} e^{-n/(k+1)}
{1\over 2(k+1)} =\sum_{k=1}^{n-1} e^{-n/(k+1)}  {n\over
2(k+1)^2}{(k+1)\over n}.
\leqno(2.49)
$$
We show that this expression converges to
$${1\over 2} \int_1^\infty e^{-\a}{1\over \a} d\a .$$

Consider the following subdivision:
$$ t_1={n\over 2},\,  \ldots \ ,t_k={n\over k+1}, \, \ldots \  ,t_{n-1}=1.$$
We have $t_{k-1}- t_k ={n\over k}-{n\over k+1}={n\over (k+1)k}$. Fix a
large integer $M$ and write
$$\eqalign{{1\over 2}\int_1^{t_M}e^{-\a}{1\over \a}
d\a&=\sum_{k=M+1}^{n-1}\int_{t_k}^{t_{k-1}}e^{-\a}{1\over 2\a}  d\a \le
\sum_{k=M+1}^{n-1} e^{-t_k}{1\over t_k}{(t_k-t_{k-1})\over 2 }\cr &
=\sum_{k=M+1}^{n-1}
e^{-{n\over k+1}}{k+1\over n}\cdot {n\over 2 (k+1)(k+1)}\cdot{k+1\over
k}\cr &\le {M+1\over M}\sum_{k=1}^{n-1}
e^{-{n\over k+1}}{n\over 2 (k+1)^2}{k+1\over n}.\cr}$$
Since $t_M={n\over M+1}\to \infty$, when $n$ tends to infinity, $M$ fixed,
we obtain
$$\liminf_{n\to \infty} \sum_{k=1}^{n-1}   {Y_k\over 2(k+1)}\ge {M\over
M+1}\cdot {1\over 2}\int_1^\infty  e^{-\a}{1\over \a} d\a.$$
By letting $M$ tend to infinity, we establish
$$\liminf_{n\to \infty} \sum_{k=1}^{n-1}   {Y_k\over 2(k+1)}\ge   {1\over
2}\int_1^\infty  e^{-\a}{1\over \a} d\a.$$
The upper bound comes similarly and we obtain
$$  \sum_{k=1}^{n-1}   {Y_k\over 2(k+1)}\ \buildrel{n\to \infty}\over
{\longrightarrow}\
{1\over 2}\int_1^\infty  e^{-\a}{1\over \a}d\a.
\leqno(2.50)$$

\smallskip
Now we turn to the last term in (2.48) showing that
$$\sum_{k=1}^{n-1}   {Z_k\over  (k+1)}= \sum_{k=1}^{n-1}  {1\over
(k+1)}\int_0^{n/(k+1)} \Big(\big(1-{v\over
n}\big)^{n-3/2} -e^{-v}\big)  dv \ \buildrel{n\to \infty}\over
{\longrightarrow}\   0.
\leqno(2.51)$$

\noi It would then follow from (2.46)--(2.51) that
$$(n-3/2)\Big[\sum_{k=1}^{n-1} (k+1)^{-(n+1/2)}B_k-\sum_{k=1}^{n-1}{1\over
2(k+1)}\Big] \ \buildrel{n\to \infty}\over {\longrightarrow}\
e^{-1}-1-{1\over 2}\int_1^\infty  e^{-\a}{1\over \a} d\a.
\leqno(2.52)$$

Let us compare the latter with the expression suggested in (2.26):
$$ \int_1^\infty \phi_2(d\a)d\a = \int_1^\infty {2e^{-\a}+\a e^{-\a}
-2\over 2\a^2}d\a=\int_1^\infty { e^{-\a} \over
 \a^2}d\a +{1\over 2}\int_1^\infty {2e^{-\a}+\a e^{-\a} -2\over 2\a^2}d\a.$$
Integration by parts yields
$$\int_1^\infty { e^{-\a} \over
 \a^2}d\a=\int_1^\infty e^{-\a}d\big({ -1  \over
 \a }\big)=e^{-\a} \big({ -1  \over
 \a }\big)\Big|_1^\infty - \int_1^\infty   \big({ -1  \over
 \a }\big)d(e^{-\a})=e^{-1}-\int_1^\infty  e^{-\a}{1\over \a}
d\a.$$
Hence,
$$ \int_1^\infty \phi_2(d\a)d\a=e^{-1}-1-{1\over 2}\int_1^\infty
e^{-\a}{1\over \a}
d\a, $$
as stated in (2.52). It remains to prove (2.51)
\smallskip
Write the sum as one integral. We have
$$\eqalign{ \bigg| \sum_{k=1}^{n-1}  {1\over  (k+1)}\int_0^{n/(k+1)}
\Big(\big(1-{v\over
n}\big)^{n-3/2}& -e^{-v}\Big)  dv \bigg|\cr &\le     \int_0^{n/(k+1)}
\Big|\big(1-{v\over
n}\big)^{n-3/2} -e^{-v}  \Big| dv \cdot\bigg(\sum_{k=1}^{n-1}  {1\over
(k+1)}\bigg)\cr&\le (\log n)
\int_0^{n/2}
\Big|\big(1-{v\over n}\big)^{n-3/2} -e^{-v}  \Big| dv\cr}.
\leqno(2.53)$$
Split the integration domain $[0,n/2]$ in $[0,A]$ and $]A,n/2]$ with
$A=A(n)$ specified below. For the second domain, we use the
elementary estimate
$$\big(1-{v\over n}\big)^{n-3/2}\le \exp\big\{-{(n-3/2)\over n}v\big\}\le
e^{-v/2}, \qq n\ge 3. $$
We thus get the estimate
$$ (\log n)({n\over 2})e^{-A/2}.$$
For the first domain, we have
$$\Big| e^{-v}- \big(1-{v\over
n}\big)^{n-3/2}    \Big|=\Big| e^{-v}- \big(1-{v\over
n}\big)^{n } \big(1-{v\over
n}\big)^{ -3/2}    \Big|:=\Big| e^{-v}- \big(1-{v\over
n}\big)^{n } h    \Big|,$$
while
$$\eqalign{\max\Big\{ e^{-v}- \big(1-{v\over
n}\big)^{n } h   ,\big(1-{v\over
n}\big)^{n } h - e^{-v} \Big\}
&\le   \max \Big\{ e^{-v}- \big(1-{v\over n}\big)^{n }, (h-1)  e^{-v} \Big\}
\cr
&\le    e^{-v}- \big(1-{v\over n}\big)^{n }  +(h-1)  e^{-v}.\cr}$$
We use the following estimate ([Mi] p.266)
$$ e^{-v}-(1-{v\over n})^n\le {v^2\over 2n}, \qq v\le n$$
and also $h=\big(1-{v\over
n}\big)^{ -3/2} \le \big(1-{A\over
n}\big)^{ -3/2} $. It follows that the expression in (2.53) is bounded by
$$(\log n)\Big({n\over 2} e^{-A/2}+ {A^3\over 6n}+ \big(1-{A\over
n}\big)^{ -3/2} -1\Big). $$
By letting $A=n^{1/4}$, we get an expression tending to $0$ when $n$ tends
to infinity. Hence we are also done with (2.26).

\vskip 0,7cm\noi \centerline{\gum 3. Final proofs}
\medskip\medskip
\medskip

\noi{\gum 3.1. Control of approximation of zeta function}
\medskip
Recall according to the notation (1.21), with $\sigma=1/2$ here,  that
$$
  Z_n(x)  = \sum_{k\le x} {1\over k^{{1\over 2}+iS_n}}
-    {x^{1-({1\over 2}+iS_n)}\over 1-({1\over 2}+iS_n)},
 $$
and put
$$\zeta_n= \zeta({1\over 2}+iS_n).
$$
We will show now that $Z_n(x)$ provides a good approximation to $\zeta_n$
in the following sense.

\medskip \noi {\gem Proposition 2.}  {\it For each positive integer $n$,}
$$ \E\big| Z_n(x)-\zeta_n\big|^2 \ \buildrel{x\to \infty}\over
{\longrightarrow}\   0.$$
    \medskip\noi

In order to prove this proposition, we need a series of simple technical
results.

\noi
Let $p_n(u)={n\over \pi(n^2+x^2)}$ denote the distribution density of $S_n$.
\medskip

\noi {\gem Lemma 1.} {\it Let $\a\in \R$ and $x\ge 1$. Then,
$$ \Big|\int_{|u|\ge x} e^{i\a u}p_n(u) du\Big|\le {C(n)\over |\a|x^2},$$
where  the constant $C(n)$ depends on $n$ only.}

\medskip\noi
{\it Proof.}
$$\int_{  x}^\infty e^{i\a u}p_n(u) du=\int_{  x}^\infty  p_n(u)
d\big({e^{i\a u}\over i\a}\big)=p_n(x){e^{i\a x}\over i\a}-\int_{
x}^\infty  p'_n(u) {e^{i\a u}\over i\a}du $$
We use the estimates
$$p_n(x)\le {C(n)\over x^2}, \qq p'_n(x)\le {C(n)\over x^3}.$$
Then
$$ \Big|\int_{  x}^\infty e^{i\a u}p_n(u) du \Big|\le {C(n)\over
x^2}{1\over |\a|}+\int_{
x}^\infty  {C(n)\over u^3} {du\over |\a|} \le {C(n)\over |\a|x^2}. $$
Applying this estimate to
$\int_{x}^\infty  e^{- i\a u}p_n(u) du
= \int_{-\infty  }^x  e^{i\a u}p_n(u) du$, we achieve our goal.\cqfd

\medskip

\noi {\gem Lemma 2.} For any fixed $n$, we have
$$ \int_{|u|\ge x} \Big|\sum_{m\le x}{1\over m^{{1\over
2}+iu}}\Big|^2p_n(u) du \ \buildrel{x\to \infty}\over {\longrightarrow}\
0.$$

\medskip\noi
 {\it Proof.} We write that
$$\Big|\sum_{m\le x}{1\over m^{{1\over 2}+iu}}\Big|^2=\sum_{m_1\le
x}\sum_{m_2\le x}{1\over m_1^{{1\over 2}+iu}}{1\over m_2^{{1\over
2}-iu}}=\sum_{m_1\le x}\sum_{m_2\le x}{1\over m_1^{{1/2} } m_2^{{1/ 2} }}
\left({m_2\over m_1 }\right)^{iu}$$
Thus
$$ \int_{|u|\ge x} \Big|\sum_{m\le x}{1\over m^{{1\over
2}+iu}}\Big|^2p_n(u) du =\sum_{m_1\le x}\sum_{m_2\le x}{1\over m_1^{{1/2} }
m_2^{{1/ 2} }}\int_{|u|\ge x} e^{iu\log \left({m_2\over m_1 }\right)
}p_n(u) du.$$

\noi We consider two cases: let $\b= 1/2$.
\medskip

If $|m_2-m_1|<m_1^\b$, then plainly
$$\Big|\int_{|u|\ge x} e^{iu\log \left({m_2\over m_1 }\right) }p_n(u) du
\Big|\le \int_{|u|\ge x}  p_n(u) du\le \int_{|u|\ge x}
{C(n)\over  u^2} du\le {C(n)\over  x}.$$
Therefore,
$$\eqalign{\sum_{m_1\le x\, ,\,m_2\le x\atop |m_2-m_1|<m_1^\b}{1\over (m_1
m_2)^{{1/2} }}
 \Big|\int_{|u|\ge x} e^{iu\log \left({m_2\over m_1 }\right) }p_n(u)
du\Big|&\le {C(n)\over  x}\sum_{m_1\le x }\sum_{ m_2\le x\atop
|m_2-m_1|<m_1^\b}{1\over (m_1 m_2)^{{1/2} }} \cr & \le {C(n)\over
x}\sum_{m_1\le x } {(2m_1^\b)\over (m_1 )^{ 1/2  }(m_1-m_1^\b)^{1/2}
}
\cr
&\le C\ {C(n)\over  x}\sum_{m_1\le x }
 m_1^{\b-1} \le {C\cdot C(n)\over  x} x^\b
 \cr
 &=C\cdot C(n)x^{\b-1}
 \cr
 &=C\cdot C(n)x^{-1/2} \ \buildrel{x\to \infty}\over
{\longrightarrow}\   0.\cr}
$$
\medskip

If $|m_2-m_1|\ge m_1^\b$,   either $m_2 - m_1\ge m_1^\b $,    then by
letting $ \psi := \log \big({m_2\over m_1 }\big)$  we get
$$|\psi|\ge \log \big({m_1+m_1^\b\over m_1 }\big)
= \log \big(1+m_1^{\b-1} \big)
\ge C m_1^{\b-1}.
$$
Or $m_1 - m_2\ge m_1^\b $, which implies $m_1 - m_2\ge m_1.m_1^{\b -1}\ge
m_2.m_1^{\b -1}$. And so
$$|\psi|= \log \big({m_1\over m_2 }\big)\ge  \log \big(1+m_1^{\b-1} \big)
\ge C m_1^{\b-1}.
$$ By applying Lemma 1, we find that
$$\Big| \int_{|u|\ge x} e^{iu\log \left({m_2\over m_1 }\right) }p_n(u)
du\Big|=\Big| \int_{|u|\ge x} e^{iu\psi }p_n(u) du\Big|\le
{C(n)\over |\psi|x^2}\le
{C(n)\over m_1^{\b-1}x^2}={C(n) \over  x^2}m_1^{1-\b},
$$
and we get
$$\eqalign{
\sum_{ m_1\le x\, ,\,m_2\le x\atop |m_2-m_1|\ge m_1^\b}
  {1\over (m_1 m_2)^{{1/2} }}
 \Big|\int_{|u|\ge x} e^{iu\log \left({m_2\over m_1 }\right) }p_n(u) du\Big|
 &\le {C(n) \over  x^2}  \sum_{ m_1\le x\, ,\,m_2\le x }
 { m_1^{1-\b }\over (m_1 m_2)^{{1/2} }}
  \cr
  &\le {C(n) \over  x^2} \sum_{ m_1\le x\, ,\,m_2\le x}  m_1^{ -\b+1/
2}m_2^{ - 1/ 2}
 \cr
 & \le {C(n) \over  x^2}\Big(\sum_{m_1 \le x }  m_1^{ -\b+1/
2}\Big)\Big(\sum_{m_2 \le x }m_2^{ - 1/ 2}
  \Big)
 \cr
    &\le {C(n) \over  x^2}   x^{ -\b+3/ 2}x^{  1/ 2}= C(n)    x^{ -\b }
 \cr
 & =C(n)    x^{ -1/2 } \ \buildrel{x\to \infty}\over {\longrightarrow}\   0.
 \cr
 }
 $$
\cqfd
\medskip

\noi {\gem Lemma 3.} \ For any fixed $n$, we have
$$ \int_{|u|\ge x} \Big| {x^{{1\over 2}-iu}\over {1\over 2}-iu
}\Big|^2p_n(u) du \
\buildrel{x\to \infty}\over {\longrightarrow}\   0.$$

 \noi
{\it Proof.}
$$\eqalign{\int_{|u|\ge x} \Big| {x^{{1\over 2}-iu}\over {1\over 2}-iu
}\Big|^2p_n(u) du& \le x\int_{|u|\ge x}{1\over |u|^2}p_n(u) du\le
C(n)x \int_{|u|\ge x}{du\over |u|^4}\le C(n)x\cdot x^{-3}\cr&={C(n)\over
x^{ 2}}\ \buildrel{x\to \infty}\over {\longrightarrow}\
0.\cr} $$
\cqfd

\noi {\it Proof of Proposition 2.} We set
$$h(x,u)=\sum_{k\le x} {1\over k^{{1\over 2}+iu}}
-    {x^{1-({1\over 2}+iu)}\over 1-({1\over 2}+iu)}. $$
Then
$$ Z_n(x)=h(x,S_n) ,$$
and we have
$$\E\big| Z_n(x)-\zeta_n\big|^2=\E\big| h(x,S_n)-\zeta({1\over
2}+iS_n)\big|^2=\int_{-\infty}^\infty \big| h(x,u)-\zeta({1\over
2}+iu)\big|^2p_n(u) du\qq\qq\qq\quad$$
$$\quad\le \int_{|u|\le x } \big| h(x,u)-\zeta({1\over
2}+iu)\big|^2p_n(u) du+ 2\int_{|u|> x } \big| h(x,u) \big|^2p_n(u)
du+2\int_{|u|> x } \big|  \zeta({1\over
2}+iu)\big|^2p_n(u) du.
$$
Concerning the first integral, we have by (1.20)
$$ \int_{|u|\le x } \big| h(x,u)-\zeta({1\over
2}+iu)\big|^2p_n(u) du\le \max_{|u|\le x}\big| h(x,u)-\zeta({1\over
2}+iu)\big|^2\le {C\over |x| }\ \buildrel{x\to \infty}\over
{\longrightarrow}\   0.
$$
The second integral,
$$
\eqalign{\int_{|u|> x } \big| h(x,u) \big|^2p_n(u) du&\le 2\int_{|u|> x }
\big|\sum_{k\le x} {1\over k^{{1\over 2}+it}}
  \big|^2p_n(u) du+2\int_{|u|> x } \big|    {x^{1-({1\over 2}+iu)}\over
1-({1\over
2}+iu)}\big|^2p_n(u) du \cr
\cr}
$$
tends to zero, as a consequence of Lemmas 2 and 3.
For controlling the third integral, we use that (see (1.9))
$$\int_{|u|\le T}|\zeta({1\over 2}+iu)|^4du \le CT(\log T)^4.
$$
We have
$$\eqalign{\int_{|u|> x } \big|  \zeta({1\over 2}+iu)\big|^2p_n(u) du& \le
\sum_{m:2^m\ge x}\int_{|u|\in [2^{m-1}, 2^m ] }
\big|  \zeta({1\over 2}+iu)\big|^2p_n(u) du\cr &\le  \sum_{m:2^m\ge x}
\Big(\max_{|u|\ge 2^{m-1}}p_n(u) \Big)\cdot\int_{|u|\in [2^{m-1}, 2^m ] }
\big|
\zeta({1\over 2}+iu)\big|^2du
\cr
&\le  \sum_{m:2^m\ge x} {C(n)\over (2^m)^2}\cdot\left(\int_{|u|\le 2^m   }
\big|  \zeta({1\over 2}+iu)\big|^4du\right)^{1/2}\cdot (2^m)^{1/2}
\cr
&\le  \sum_{m \ge {\log x\over \log 2}}    {C(n)\over (2^m)^2}\,  \cdot \,
\Big(2^m [m\, \log 2]^4 \Big)^{1/2} \cdot
(2^m)^{1/2}
\cr
&\le  C \ \cdot \ C(n) \sum_{m \ge {\log x\over \log 2}} { m^2\over  2^m
}\ \buildrel{x\to \infty}\over
{\longrightarrow}\   0,
\cr}
$$
and the proof is now complete.\cqfd

\noi Note that a weaker result than (1.9), for instance Theorem 7.4 in [T]
asserting that $\int_0^T |\zeta({1\over 2}+ it)|^2{\rm d}t\sim T\log T$ as
$T$
tends to infinity, would have also be suitable.
\medskip

\noi{\gum 3.2. Proof of Theorem 1.}
\medskip

\noi  In this step, we prove Theorem 1 by giving the estimates of the
covariance of the system $\{ {\cal Z}_n,n\ge 1\}$.
Recall that
$${\cal Z}_n = {\zeta}({1/  2}+iS_n) -\E\,  {\zeta}({1/ 2}+iS_n)=\zeta_n-\E
\zeta_n. $$
We approximate $\zeta_n$ by $Z_n(x)$. By Proposition 2 we know that
$$ \E\, \big| Z_n(x)-\zeta_n\big|^2 \ \buildrel{x\to \infty}\over
{\longrightarrow}\   0.   \leqno(3.1)
$$
On the other hand, by using (2.21), (2.14), (2.16), we obtain from (2.3)
$$\eqalign{
\E\, |Z_n(x)|^2 &=  {2x\over n+1/2} + K_n + {2x\over n+1/2}
 -  {4x\over n+1/2} + {2\over 2n-1} + o(1)
\cr
& = K_n  + {1\over n-1/2} + o(1), \qquad x\to\infty,
\cr},
$$
where $K_n$ is given in (2.21). Hence,
$$\E\,   |\zeta _n|^2=K_n  + {1\over n-1/2}<\infty.
$$
It follows from (3.1) that
$$\E\,   \zeta _n \overline{   \zeta} _m =\lim_{x\to \infty}   \E\, Z_n(x)
\overline{  Z}_m(x) . \leqno(3.2)
$$
It also follows from (3.1) and (2.4) that
$$\E\,   \zeta_n   =\lim_{x\to \infty}   \E\, Z_n(x)   =\zeta({1\over 2}+n)
-   {8n\over 4n^2-1}.
\leqno(3.3)
$$
Since
$$ \E\, {\cal Z}_n\overline{{\cal Z}_m}=
\E\,   \zeta _n \overline{   \zeta} _m  -\E\,   \zeta_n \overline{ \E\,
\zeta _m},
$$
we obtain from (3.2) and (3.3)
$$ \E\, {\cal Z}_n\overline{{\cal Z}_m}
= \lim_{x\to \infty}   \E\, Z_n(x)  \overline{  Z}_m(x)-\big(\zeta({1\over
2}+n) -   {8n\over 4n^2-1}\big)\big(\zeta({1\over 2}+m) -   {8m\over
4m^2-1}\big)
. \leqno(3.4)
$$
In particular,
$$\E\,   | {\cal Z}_n|^2=K_n  + {1\over n-1/2}-\big(\zeta({1\over 2}+n) -
{8n\over 4n^2-1}\big)^2,  \leqno(3.5)
$$
and the first claim of Theorem 1 follows.
\smallskip
By (2.3),
$
\E\,  Z_n(x) \bar Z_m(x)
  =\E\,  Z_{n1} \bar Z_{m1} -\E\, Z_{n1} \bar Z_{m2} -\E\,  Z_{n2} \bar
Z_{m1} +\E\,  Z_{n2} \bar Z_{m2}$.
Recall that we proved in (2.20), (2.15), (2.17) and (2.13) respectively,
for $m>n+1$, as $x$ tends to infinity, that
$$\eqalign{  \E\, Z_{n1}\bar
Z_{m1}
 &= \zeta((m-n)+1)+\theta \left(  {1\over m-1/2} +  {1\over n-1/2}
\right)\zeta(m-n) +o(1),
\cr
\E\, Z_{n1}\bar Z_{m2}&=
{-2(m-n)\zeta(n+1/2)\over (m+1/2)(2n-m+1/2)}
 + \
{2n \ \zeta(m-n)\over (m-1/2)(2n-m+1/2)}
+ o(1),\cr \E\,  Z_{n2} \bar
Z_{m1}& =\overline{\E\, Z_{m1}\bar Z_{n2}} = {2n \ \zeta(m-n)\over n^2-1/4}
+ o(1), \cr \E\, Z_{n2}\bar Z_{m2}& =
 {4n(m-n)\over ((m-n)^2-1)(n^2-1/4)}+o(1),\cr}$$
where $\theta=\theta(n,m)\in [0,1]$.
We get in view of  (2.4)
$$\eqalign{\Big|\E\, Z_n(x)  \overline{  Z}_m(x)-\E\, Z_n(x)  \overline{\E\,
Z_m(x)}\Big|&\le \Big|\zeta((m-n)+1) \cr &\ -\big(\zeta({1\over 2}+n) -
 {8n\over 4n^2-1}\big)\big(\zeta({1\over 2}+m) -   {8m\over
4m^2-1}\big)\Big|\cr &\
+
\left(  {1\over m-1/2} +  {1\over n-1/2}
\right)\zeta(m-n)
\cr &\  +\bigg|{ 2(m-n)\zeta(n+1/2)\over (m+1/2)(2n-m+1/2)}-
{2n \ \zeta(m-n)\over (m-1/2)(2n-m+1/2)}
\cr &\quad  -{2n \ \zeta(m-n)\over n^2-1/4}+{4n(m-n)\over
((m-n)^2-1)(n^2-1/4)}  \bigg|
+o(1).\cr}
$$
By using (3.4) and letting $x$ tend to infinity,
we obtain for any fixed pair of integers $n,m$ with $m>n+1$
$$\eqalign{\Big|\E\,{\cal Z}_n\overline{{\cal Z}_m}\Big|&\le
\Big|\zeta((m-n)+1) -\big(\zeta({1\over 2}+n) -
 {8n\over 4n^2-1}\big)\big(\zeta({1\over 2}+m) -   {8m\over 4m^2-1}\big)\Big|
\cr & \quad +  \left(  {1\over m-1/2} +  {1\over n-1/2}
\right)\zeta(m-n)
\cr &\quad  +\bigg|{ 2(m-n)\zeta(n+1/2)\over (m+1/2)(2n-m+1/2)}-
{2n \ \zeta(m-n)\over (m-1/2)(2n-m+1/2)}
\cr &\quad  -{2n \ \zeta(m-n)\over n^2-1/4}+{4n(m-n)\over
((m-n)^2-1)(n^2-1/4)}  \bigg|
 .\cr}$$

\noindent But
$$\zeta((m-n)+1)  -\big(\zeta({1\over 2}+n) -
 {8n\over 4n^2-1}\big)\big(\zeta({1\over 2}+m) -   {8m\over 4m^2-1}\big)
=\zeta((m-n)+1) $$
$$-\zeta({1\over 2}+n)\zeta({1\over 2}+m) +{8m\over 4m^2-1}\cdot
\zeta({1\over 2}+n)+
{8n\over 4n^2-1}\cdot \zeta({1\over 2}+m) -{64nm\over (4n^2-1)(4m^2-1)}.$$
Now $$\zeta((m-n)+1)-\zeta({1\over 2}+n)\zeta({1\over 2}+m)
  =\sum_{k=1}^\infty {1\over
k^{(m-n)+1}}-\sum_{k=1}^\infty\sum_{\ell=1}^\infty {1\over k^{{1\over
2}+n}\ell^{{1\over 2}+m}}
$$
$$ =\sum_{k=2}^\infty {1\over k^{(m-n)+1}}
  -\sum_{k=2}^\infty{1\over k^{{1\over 2}+n} }
  -\sum_{\ell=2}^\infty {1\over \ell^{{1\over 2}+m}}
  -\sum_{k=2}^\infty{1\over k^{{1\over 2}+n} }
   \sum_{\ell=2}^\infty {1\over  \ell^{{1\over 2}+m}}\ .
$$
And since for any $D>1$,
$$\sum_{k=2}^\infty {1\over k^{D}}\le  {1\over 2^{D}}\Big(1+{2\over
D-1}\Big),$$
it follows that
$$\Big|\zeta((m-n)+1)  -\zeta({1\over
2}+n)\zeta({1\over 2}+m)\Big|\le \! C\max\Big({1\over 2^{(m-n)+1}},{1\over
2^{{1\over 2}+n}},{1\over 2^{{1\over 2}+m}}\Big)\le
C\max\Big({1\over 2^{ m-n  }},{1\over 2^{ n}} \Big).$$
Further
$$\sup_{m>n+1}  { m\over 4m^2-1}  ={\cal O}({1\over n}), \qq { n\over
4n^2-1}  ={\cal O}({1\over n}), $$
so that
$$\sup_{m>n+1}\bigg|{8m\over 4m^2-1}\cdot \zeta({1\over 2}+n)+
{8n\over 4n^2-1}\cdot \zeta({1\over 2}+m) -{64nm\over
(4n^2-1)(4m^2-1)}\bigg| ={\cal O}({1\over n}). $$
 For the other terms, we have uniformly over $m$ such that $m>n+1$
 $$\eqalign{  \left(  {1\over m-1/2} +  {1\over n-1/2}
\right)\zeta(m-n) &={\cal O} ( {1\over n}),
\cr
    {2n \ \zeta(m-n)\over n^2-1/4}
&={\cal O} ( {1\over n}) ,
\cr
 {4n(m-n)\over ((m-n)^2-1)(n^2-1/4)}&={\cal O} ( {1\over n}) .\qq\qq\cr}$$
Consider finally the last term
$$ {-2(m-n)\zeta(n+1/2)\over (m+1/2)(2n-m+1/2)}
 +
{2n \, \zeta(m-n)\over (m-1/2)(2n-m+1/2)}
 $$
We have
$${2n \over (m-1/2)(2n-m+1/2)}
-{2(m-n) \over (m-1/2)(2n-m+1/2)}=2{(2n-m)\over (m-1/2)(2n-m+1/2)} $$
thus
$$\eqalign{\Big|{2n \over (m-1/2)(2n-m+1/2)}
-{2(m-n) \over (m-1/2)(2n-m+1/2)}\Big|&\le {2\over  m-1/2 }\max_{U=2n-m\in
\Z}{|U|\over   |U+1/2| }\cr &\le C/m. \cr}$$
Further
$${2(m-n) \over (m-1/2)(2n-m+1/2)}-{2(m-n) \over (m+1/2)(2n-m+1/2)}={1
\over (m^2-1/4) }.{2(m-n) \over
 (2n-m+1/2)} .$$ Now observe that the function $f(A):= {A\over n-A+1/2} $
defined for $A$ integer,
 has maximal absolute value less than $Cn$. Hence, as
$f(m-n)={ (m-n) \over (2n-m+1/2)}$, we deduce
$$\max_{m>n+1}\Big|{2(m-n) \over (m-1/2)(2n-m+1/2)}-{2(m-n) \over
(m+1/2)(2n-m+1/2)}\Big| \le C\max_{m>n+1}{n \over m^2}\le
  {C \over n}\ .
$$
Consequently,
$$\max_{m>n+1}\Big|{2n \over (m-1/2)(2n-m+1/2)}
-{2(m-n) \over (m+1/2)(2n-m+1/2)}\Big| \le
  {C \over n}\ .
$$
Write that
$$\displaylines{ {-2(m-n)\zeta(n+1/2)\over (m+1/2)(2n-m+1/2)}
 +
{2n \, \zeta(m-n)\over (m-1/2)(2n-m+1/2)}
\hfill\cr\hfill ={ 2(m-n) \over (m+1/2)(2n-m+1/2)}  \big\{\zeta(m-n)
-\zeta(n+1/2)\big\}
 \cr\hfill +
\zeta(n+1/2)\Big\{{2n \,  \over (m-1/2)(2n-m+1/2)} -{2(m-n) \over
(m+1/2)(2n-m+1/2)}\Big\}.\cr}$$
We already know that the absolute value of the last term is less than $ {C
/ n}.$ To control the first term we proceed as
before: since
$$\big|\zeta(m-n) -\zeta(n+1/2)\big|\le C  \max \Big({1\over 2^{ m-n
}},{1\over 2^{n  }}   \Big)$$
and
$$
\big|{ 2(m-n) \over (m+1/2)(2n-m+1/2)}\big|={2 |f(m-n)|\over  (m+1/2)}\le
C{n\over m} \le C,
$$
we get
$$\Big|{ 2(m-n) \over (m+1/2)(2n-m+1/2)} \Big\{\zeta(m-n)
-\zeta(n+1/2)\Big\}\Big|
\le  C
\, \max \Big({1\over 2^{ m-n
}},{1\over 2^{n  }}   \Big)
.
$$
Therefore, for $m>n+1$
$$ \big|\E\,{\cal Z}_n\overline{{\cal Z}_m}\big|\le C \max \Big( {1\over
n},{1\over 2^{ m-n  }}   \Big),$$
as claimed in Theorem 1.
\cqfd
\medskip

\noi {\gum 3.3. Asymptotic behavior along the Cauchy random walk}
\medskip

\noi
In this subsection, we give the proof of Theorem 2.
The essential step consists of controlling the increments
$$\E\, \Big|\sum_{i\le n\le j\atop n\, even}{\cal Z}_n\Big|^2, \qq  \E\,
\Big|\sum_{i\le n\le j\atop n\, odd}{\cal Z}_n\Big|^2. $$
Since the two increments are treated in exactly the same way, we only
consider the first one.
We use Theorem 1. By developing the sum
$$\E\, \Big|\sum_{i\le n\le j\atop n\, even}{\cal Z}_n\Big|^2=   \sum_{i\le
n\le j\atop n\, even}\E\,|{\cal Z}_n |^2+
 \sum_{i\le n\le j\atop n\, even}\sum_{i\le m\le j\atop m\, even}\E\,{\cal
Z}_n\overline{{\cal Z}_m}\le C \sum_{i\le n\le j\atop n\,
even}\log n  +
C \!\sum_{i\le n<m\le j\atop n,m\, even}  \max \big( {1\over n},{1\over 2^{
m-n  }}   \big).$$
But
$$\sum_{i\le n<m\le j\atop n,m\, even}   {1\over n} \le \Big(\sum_{  n \le
j }{1\over n}\Big)\Big(\sum_{i\le   m\le
j\ }  1\Big)\le C(\log j)(j-i)
$$
and
$$\sum_{i\le n<m\le j\atop n,m\, even}   2^{-( m-n)  }\le \Big(\sum_{i\le n
\le j }1 \Big)\Big(\sum_{  m>n} 2^{-( m-n)
}\Big)\le  (j-i)\Big(\sum_{  h\ge 1} 2^{-h
}\Big)= C (j-i).
$$
Therefore,
$$\E\, \Big|\sum_{i\le n\le j\atop n\, even}{\cal Z}_n\Big|^2\le C(\log
j)(j-i).$$
And by operating the same way for the odd part, we get that there exists a
constant $C$ such that for any $j>i$,
$$\E\, \Big|\sum_{i\le n\le j }{\cal Z}_n\Big|^2\le C(\log j)(j-i).$$
Now the conclusion of Theorem 2 is easily obtained from Theorem 1.10 in
[W2], which we recall now.
\medskip

\noi {\bf Proposition 3.}\ {\it  Let $\{m_l, l\ge 1\}$ be a  sequence of
positive reals with partial sums $M_n=\sum_{l=1}^nm_l$
tending to infinity with $n$. Assume that
$$\log{M_n\over m_n}\sim \log M_n.$$
  Let $\Phi :\R^+\rightarrow\R^+$ be a concave nondecreasing function. Then
any sequence
  $\{\xi_l,l\ge 1\}$ of random variables satisfying the  increment condition
$$\E\big| \sum_{l=i}^j \xi_l\big|^2\le \Phi(\sum_{l=1}^j m_l) \,
\big(\sum_{l=i}^j m_l\big), \qq \qq (i\le j) $$
  also verifies  for any $\tau>3/2$}
$$ {  \sum_{l=1}^n\xi_l \over \Phi(M_n )^{1/2}\log^{ \tau }(1+M_n)}\
\buildrel{a.s.}\over{\longrightarrow}\  0\qquad{\it and
}\qquad
\Big\| \,\sup_{n\ge 1} {\big| \sum_{l=1}^n\xi_l\big|\over \Phi(M_n
)^{1/2}\log^{ \tau }(1+M_n)}\, \Big\|_2<\infty.  $$
We apply this result to ${\cal Z}_n$ with the choice  $m_l\equiv 1$ and
$\Phi(x)=\log(1+x)$ and obtain the assertion
of Theorem 2.
\cqfd
\bigskip

{\bf Acknowlegements.} \ We a grateful to Prof. E. Bombieri for informing
us about actual references concerning the
introductory part.  The work of the first named author was supported by
grants NSh.422.206.1, RFBR 05-01-00911,
and INTAS 03-51-5018.
\bigskip

\medskip\noi
   {\gum References}
\smallskip

\noi  [Bi]  {Biane Ph.} [2003]  {\sl La fonction z\^eta de Riemann et les
probabilit\'es}.
In: La fonction z\^eta. Ed. \'Ec.Polytech., 165--193.

\smallskip\noi [Bl]  Blanchard A. [1969] {\sl Initiation
\`a la th\'eorie analytique des nombres premiers},  Travaux et Recherches
Math\'ematiques,
{\bf 19}, Dunod, Paris.

\smallskip\noi  [Bo]  {Bohr H.} [1952]  {\sl Collected Mathematical Works}.
Dansk Mat. Forening,
Copenhagen.

\smallskip\noi  [BJ]  {Bohr H., Jessen B.} [1930/1932]  {\sl \"Uber die
Wertverteilung der
Riemannsche Zetafunktion}, Erste Mitteilung. Acta Math. {\bf 54},
1--35; Zweite Mitteilung, ibid. {\bf 58}, 1--55.


\smallskip\noi  [BPY]  {Biane Ph., Pitman J., Yor M.} [2001]  {\sl
Probability laws related
to the Jacobi theta and Riemann zeta functions, and Brownian excursions}.
Bull. Amer. Math. Soc. {\bf 38}, no 4, 435--465.

 \smallskip\noi [Gh]  Ghosh A. [1983] {\sl On the Riemann zeta-function --
mean value theorems
and the distribution of $|S(T)|$}. J. Number Theory {\bf 17}, 93--102.

\smallskip\noi  [GH]  {Gelbart S., Miller S.} [2004]  {\sl Riemann's zeta
function and beyond}.
Bull. Amer. Math. Soc. {\bf 41}, 59--112.

\smallskip\noi [H1]  Huxley M. N. [1972] {\sl The Distribution of Prime
Numbers}. Oxford Math. Monographs.

\smallskip\noi  [H2]  {Huxley, N. N. } [2005]  {\sl
Exponential sums and the Riemann zeta function},  V. Proc. London Math.
Soc. (3) {\bf 90}  no. 1, 1--41.

\smallskip\noi [HM]  Hattori T., Matsumoto K. [1999] {\sl A limit theorem for
Bohr-Jessen's probability measures of the Riemann zeta-function}.
J. Reine Angew. Math. {\bf 507}, 219--232.

\smallskip\noi [I1] Ivi\'c  A. [2003] {\sl The Riemann Zeta-function,
Theory and Applications.
} Dover Publications, Inc., Mineola, NY.

\smallskip\noi [I2] Ivi\'c  A. [2005] {\sl The modified Mellin transform of
powers
of the zeta-function}. Preprint.

\smallskip\noi [Jo] Joyner D.  [1986] {\sl Distribution Theorems of
$L$-functions},
Longman Scientific \& Technical.


\smallskip\noi [La]  Laurin\v cikas A. [2002] {\sl A probabilistic equivalent
of the Lindel\" of Hypothesis}. In: Analytic and Probabilistic Methods in
Number Theory.
Palanga, 157--161.

\smallskip\noi [LS] Laurin\v cikas A., Steuding J.  [2003]
{\sl A short note on the Lindel\"of Hypothesis}. Lithuanian Math. J.
{\bf 43}, no 1, 51--55.

 \smallskip\noi [Mi] Mitrinovi\'c D. S. [1970] {\sl Analytic Inequalities.}
 Springer Verlag, Berlin-New York.

\smallskip\noi [Mo] Montgomery H. L. [1973] {\sl The pair correlation of zeros
of the zeta function}, Proc. Sympos. Pure Math. {\bf XXIV}, Analytic Number
Theory,  181--193.

\smallskip\noi [Po] P\'olya G.  [1974] {\sl Collected Papers. Vol II:
Location of Zeros}.
Mathematicians of Our Time, {\bf 8}, The MIT Press, Cambridge,Mass.--London.

\smallskip \noi [T] {Titchmarsh E. C.} [1951]   {\sl The Theory of the
Riemann-Zeta Function},
Oxford University Press.

\smallskip \noi [T] {Vinogradov I. M.} [1958]   {\sl   A new estimate of the
function $\zeta (1+it)$},  (Russian) Izv. Akad. Nauk SSSR. Ser. Mat. {\bf 22},
 161--164.
\smallskip\noi [W1]  Weber M.,  [2005]  {\sl Divisors, spin sums  and
the functional equation of the Zeta-Riemann function}, Periodica Math.
Hungar.
{\bf 51}, no 1, 119--131.

\smallskip\noi [W2] Weber M.  [2006] {\sl Uniform bounds under increment
conditions}.
Trans. Amer. Math. Soc. {\bf 358}, no.2, 911--936.
\vskip 15pt

\bigskip

\noi {\phh Mikhail A. Lifshits
\par\noi  Department of Mathematics and Mechanics,
\par\noi St.Petersburg State University,
\par\noi 198504, Stary Peterhof, Russia.
\noi E-mail:\ \tt lifts@mail.rcom.ru}
\medskip

\noi {\phh Michel J. G. Weber \par\noi  Math\'ematique (IRMA),
Universit\'e Louis-Pasteur et C.N.R.S.,   7  rue Ren\'e Descartes,
\par \noi  67084 Strasbourg Cedex, France. E-mail: \  \tt
weber@math.u-strasbg.fr}
 \bye